\documentclass[12pt,oneside]{article}

\usepackage{amsmath}
\usepackage{amssymb}
\usepackage{color}
\usepackage{rotating}

\usepackage{graphicx}
\usepackage{subfigure}

\usepackage{url}

\usepackage{appendix}

\usepackage{longtable}

\author{
Zhiqi Liu\footnote{E-mail addresses: salhhhhh@163.com, 2020010316@buct.edu.cn,\\  2024201119@mail.buct.edu.cn}~ and 
Hui Zhou\footnote{E-mail addresses: zhouh06@qq.com, zhouhlzu06@126.com}~\footnote{Corresponding author.} 
\vspace{1em}\\ 
\footnotesize{College of Mathematics and Physics, }\\ 
\footnotesize{Beijing University of Chemical Technology, Beijing 100029, China}
}

\title{\Large Homogeneous linear recurrence relations of the determinants of distance matrices of trees}

\def\diam{{\rm diam}}
\def\K{{\sf K}}

\newtheorem{theorem}{Theorem}[section]

\newcommand*{\QEDA}{\hfill\ensuremath{\blacksquare}}  

\date{August 31, 2024}

\begin{document}

\maketitle

\begin{abstract}
In 1971, by induction on $n$ and using a two-term linear recurrence relation, Graham and Pollak got a beautiful formula \begin{center}$\det(D_n)=-(n-1)(-2)^{n-2}$\end{center} on the determinant of distance matrix $D_n$ of a tree $T_n$ on $n$ vertices. The recurrence relations are very crucial when proving this formula by inductive method: in 2006, Yan and Yeh used two-term and three-term recurrence relations; in 2020, Du and Yeh used a homogeneous linear three-term recurrence relation. In this paper, we analyze the subtree structure of the tree and find four-term, five-term, six-term and seven-term homogeneous linear recurrence relations on $\det(D_n)$, as a corollary new proofs of Graham and Pollak's formula can be given. 

\end{abstract}




\section{Introduction}\label{sec-Introduction}

In this paper, all graphs are finite, simple and undirected. The study on graphs and matrices is an important topic in algebraic graph theory. A tree is a connected graph without cycles, and it has many beautiful combinatorial properties. 

Let $n\geqslant 1$ be an integer, and let $[n]=\{1,2,\ldots,n\}$. Let $T_n$ be a tree with vertex set $[n]$. Then for any two vertices $i,j\in [n]$, there is a unique path $P_{i,j}$ connecting $i$ and $j$. The length (number of edges) of $P_{i,j}$ is called the distance between $i$ and $j$ and is denoted by $d_{i,j}$. For each $i\in [n]$, let $d_{i,i}$ be zero. Then the matrix $D_n=(d_{i,j})_{i,j\in [n]}$ is the distance matrix of $T_n$. The determinant and inverse of the distance matrix of a graph is of great interest. This kind of research is initiated by Graham and Pollak~\cite{GrahamPollak1971} on the determinant of the distance matrix of a tree, i.e., they got 
\begin{center}
  $\det(D_n)=-(n-1)(-2)^{n-2}$, 
\end{center} 
which only depends on the number $n$ of vertices of the tree and is independent with the structure of the tree. We usually call it Graham and Pollak's formula. More related results can be seen in~\cite{BapatKirklandNeumann2005,BapatLalPati2006,BapatLalPati2009,EdelbergGareyGraham1976,
GrahamHoffmanHosoya1977,GrahamLovasz1978,ZhouDing2016LAA}. New proofs of Graham and Pollak's formula can be seen in~\cite{DuYeh2020,YanYeh2006}.

We consider a decomposition of the tree $T_n$ into a subtree $S_m$ with $m~(2\leqslant m<n)$ vertices and a subtree $R$ with $n-m+1$ vertices such that $S_m$ and $R$ has exactly one common vertex $x$, and we denote it by \begin{center}$T_n=S_m\circ_x R$\end{center} or $T_n=S_m\circ R$ when the common vertex is not specified. Decompositions $T_n=S_m\circ_x R$ and $T_n=S_m\circ_y R$ are equivalent if one of the distance matrices corresponding to them can be obtained by interchanging row $x$ and row $y$ and interchanging column $x$ and column $y$ of the other distance matrix. Distance matrices corresponding to equivalent decompositions are regarded as to be of the same form, since the determinants of these distance matrices are equal. 

For a given $m$, according to different structures of $S_m$ and different common vertex $x$, we can get all the possible forms of $D_n$, and then we have the recurrence relations of $D_n$; and hence a proof of Graham and Pollak's formula would be given by induction on $n$ through these recurrence relations and initial conditions. 

The subtree $S_2$ on two vertices is a path $P_2$ with one edge $\{1,2\}$, the decompositions $T_n=P_2\circ_1 R$ and $T_n=P_2\circ_2 R$ are equivalent, so they induce the same form of $D_n$. In Graham and Pollak's proof~\cite{GrahamPollak1971}, they consider this decomposition recursively and got a two-term recurrence relation $\det(D_n)+2\det(D_{n-1})=-(-2)^{n-2}$. Combining with initial conditions $\det(D_1)=0$ and $\det(D_2)=-1$, it is not hard to solve this recurrence equation and get the Graham and Pollak's formula.  

In Du and Yeh's proof~\cite{DuYeh2020}, they consider the decomposition $T_n=S_3\circ R$. We know the subtree $S_3$ on three vertices is a path $P_3$ with edges $\{1,2\}$ and $\{2,3\}$. The decompositions $T_n=P_3\circ_1 R$ and $T_n=P_3\circ_3 R$ are equivalent, and they induce the same form of $D_n$. So there are two different forms of $D_n$, they calculated these two forms and luckily only one homogeneous linear three-term recurrence relation $\det(D_n)+4\det(D_{n-1})+4\det(D_{n-2})=0$ is derived. The next step is the same as above, combining with initial conditions $\det(D_1)=0$ and $\det(D_2)=-1$, it is not hard to solve this recurrence equation and get the Graham and Pollak's formula.

Inspired by the above proofs, we are interested in the decompositions $T_n=S_m\circ R$ and all the possible recurrence relations corresponding to different structures of the subtree $S_m$ for $m\geqslant 4$. Hence combining with suitable initial conditions $\det(D_r)$ for $r=0,1,\ldots, t_m$ where $1\leqslant t_m\leqslant m-1$, we can solve recurrence relations and get new proofs of Graham and Pollak's formula. 

In Section~\ref{sec-S4}, we consider the decomposition $T_n=S_4\circ R$. The subtree $S_4$ has two types, they yield four different forms of $D_n$ according to different common vertices, and hence we get homogeneous linear recurrence relations (including four-term recurrence relations), see in Table~\ref{tab-S4}.  
In Section~\ref{sec-S5}, we consider the decomposition $T_n=S_5\circ R$. The subtree $S_5$ has three types, they yield nine different forms of $D_n$ according to different common vertices, and hence we get homogeneous linear recurrence relations (including five-term recurrence relations), see in Table~\ref{tab-S5}. 
In Section~\ref{sec-S6}, we consider the decomposition $T_n=S_6\circ R$. The subtree $S_6$ has six types, they yield twenty different forms of $D_n$ according to different common vertices, and hence we get homogeneous linear recurrence relations (including six-term recurrence relations), see in Table~\ref{tab-S6}. 
In Section~\ref{sec-S7}, we consider the decomposition $T_n=S_7\circ R$. The subtree $S_7$ has eleven types, they yield forty-eight different forms of $D_n$ according to different common vertices, and hence we get homogeneous linear recurrence relations (including seven-term recurrence relations), see in Table~\ref{tab-S7}. 
Hence we have the following result.

\begin{theorem}\label{thm-4567-term-recurrence-relations}
For $4\leqslant m\leqslant 7$, a homogeneous linear $m$-term recurrence relation on $\det(D_n)$ can be deduced from some decomposition $T_n=S_m\circ R.$
\end{theorem}

So a question arises naturally: for which subtree $S_m$, we can deduce a homogeneous linear $m$-term recurrence relation on $\det(D_n)$ (especially for $m>7$)? Maybe the path $P_m$ is a good choice by our tables.


\section{Preliminaries}\label{sec-preliminaries}

We use $e$ to denote an appropriate size column vector whose entries are ones. Let $\alpha$ be a vector and $M$ be a matrix, we use $\alpha^T$ to denote the transpose vector of $\alpha$ and $M^T$ to denote the transpose matrix of $M$. For a matrix or a partitioned matrix $M$, we use $R_i(M)$ or $R_i$ to denote the $i$-th row or partitioned row of $M$, and $C_i(M)$ or $C_i$ to denote the $i$-th column or partitioned column of $M$; in $\det(M)$, we use $aR_i\rightarrow R_j$ to denote the elementary operation of adding $a$ times of the $i$-th row $R_i$ to the $j$-th row $R_j$, and $aC_i\rightarrow C_j$ to denote the elementary operation of adding $a$ times of the $i$-th column $C_i$ to the $j$-th column $C_j$, where $a$ is a scalar. For example, $R_3-2R_5\rightarrow R_4$ represents adding $R_3$ to $R_4$ and adding $-2$ times $R_5$ to $R_4$.

Let $m$ and $r_1,r_2,\ldots,r_m$ be positive integers, and let $P_i$ be a path of length $r_i$ and $v_i$ is an end-vertex of $P_i$ for $1\leqslant i\leqslant m$. A star-like tree $K(r_1,r_2,\ldots,r_m)$ is obtained from paths $P_1,P_2,\ldots,P_m$ by identifying the end-vertices $v_1,v_2,\ldots,v_m$. For example, a star $K_{1,3}$ is a star-like tree $K(1,1,1)$.


\section{Recurrence relations related to $T_n=S_4\circ R$}\label{sec-S4}

In this section, we consider the decomposition $T_n=S_4\circ R$ and recurrence relations related. We may suppose the vertex set of the subtree $S_4$ on four vertices is $[n]\setminus [n-4]=\{n-3,n-2,n-1,n\}$. We know the subtree $S_4$ is either a path $P_4$ or a star $K_{1,3}$. 

We may suppose $P_4$ has edges \begin{center}$\{n-3,n-2\}$, $\{n-2,n-1\}$ and $\{n-1,n\}$.\end{center} The decompositions $T_n=P_4\circ_{n} R$ and $T_n=P_4\circ_{n-3} R$ are equivalent, decompositions $T_n=P_4\circ_{n-1} R$ and $T_n=P_4\circ_{n-2} R$ are equivalent, so there are two different decompositions with form $T_n=P_4\circ R$. 

We may suppose $K_{1,3}$ has edges \begin{center}$\{n-3,n-2\}$, $\{n-3,n-1\}$ and $\{n-3,n\}$.\end{center} The decompositions $T_n=K_{1,3}\circ_{n} R$, $T_n=K_{1,3}\circ_{n-1} R$ and $T_n=K_{1,3}\circ_{n-2} R$ are equivalent, so there are two different decompositions with form $T_n=K_{1,3}\circ R$. Hence there are four different decompositions of the form $T_n=S_4\circ R$: 
\begin{center}
  $T_n=P_4\circ_{n-3} R$, $T_n=P_4\circ_{n-2} R$,\\ $T_n=K_{1,3}\circ_{n-3} R$, $T_n=K_{1,3}\circ_{n-2} R$.
\end{center} 

For the convenience of calculations in the following subsections, we usually give new labels of vertices of $S_4$ such that label $n-3$ as the common vertex of the decomposition $T_n=S_4\circ R$, and label other vertices of $S_4$ by distances from $n-3$.


\numberwithin{equation}{section}
\setcounter{equation}{0}

\subsection{The recurrence relation related to $T_n=P_4\circ_{n-3} R$}\label{subsec-S4-P4n-3}

In this subsection, we consider the decomposition $T_n=P_4\circ_{n-3} R$ of $T_n$, then the partitioned matrix $D_n$ has the following form:

\begin{equation}\label{eqn-S4-P4n-3} 
D_n =\left( \begin{matrix}
		{{D}_{n-4}} & \alpha  & \alpha +e & \alpha +2e & \alpha +3e  \\
		{{\alpha }^{T}} & 0 & 1 & 2 & 3  \\
		{{\alpha }^{T}}+{{e}^{T}} & 1 & 0 & 1 & 2  \\
		{{\alpha }^{T}}+2{{e}^{T}} & 2 & 1 & 0 & 1  \\
		{{\alpha }^{T}}+3{{e}^{T}} & 3 & 2 & 1 & 0  \\
	\end{matrix} \right),
\end{equation}
where $\alpha$ is a column vector whose entries are $d_{i,n-3}$ for $i\in [n-4]$. 

We do $R_3-2R_4\rightarrow R_5$ on $\det(D_n)$ and do $C_2-C_3\rightarrow C_4$ on the first determinant of the second line, then we get \begin{equation}\label{eqn-S4-P4n-3-4by4} 
\begin{aligned}
 \det \left( {{D}_{n}} \right)&=\left| \begin{matrix}
		{{D}_{n-4}} & \alpha  & \alpha +e & \alpha +2e & \alpha +3e  \\
		{{\alpha }^{T}} & 0 & 1 & 2 & 3  \\
		{{\alpha }^{T}}+{{e}^{T}} & 1 & 0 & 1 & 2  \\
		{{\alpha }^{T}}+2{{e}^{T}} & 2 & 1 & 0 & 1  \\
		0 & 0 & 0 & 2 & 0  \\
	\end{matrix} \right| \\ 
	& =-2\left| \begin{matrix}
		{{D}_{n-4}} & \alpha  & \alpha +e & \alpha +3e  \\
		{{\alpha }^{T}} & 0 & 1 & 3  \\
		{{\alpha }^{T}}+{{e}^{T}} & 1 & 0 & 2  \\
		{{\alpha }^{T}}+2{{e}^{T}} & 2 & 1 & 1  \\
	\end{matrix} \right| \\
 &=-2\left| \begin{matrix}
		{{D}_{n-4}} & \alpha  & \alpha +e & \alpha +2e  \\
		{{\alpha }^{T}} & 0 & 1 & 2  \\
		{{\alpha }^{T}}+{{e}^{T}} & 1 & 0 & 3  \\
		{{\alpha }^{T}}+2{{e}^{T}} & 2 & 1 & 2  \\
	\end{matrix} \right|.  
\end{aligned}
\end{equation}
Based on the property of determinants, we construct $\det \left( {{D}_{n-1}} \right)$ and expand it along the last column. This transforms the above equation into 

\begin{equation*}
\begin{aligned}
&	\det \left( {{D}_{n}} \right) =-2\left| \begin{matrix}
		{{D}_{n-4}} & \alpha  & \alpha +e & \alpha +2e  \\
		{{\alpha }^{T}} & 0 & 1 & 2  \\
		{{\alpha }^{T}}+{{e}^{T}} & 1 & 0 & 1  \\
		{{\alpha }^{T}}+2{{e}^{T}} & 2 & 1 & 0  \\
	\end{matrix} \right|-2\left| \begin{matrix}
		{{D}_{n-4}} & \alpha  & \alpha +e & 0  \\
		{{\alpha }^{T}} & 0 & 1 & 0  \\
		{{\alpha }^{T}}+{{e}^{T}} & 1 & 0 & 2  \\
		{{\alpha }^{T}}+2{{e}^{T}} & 2 & 1 & 2  \\
	\end{matrix} \right| \\ 
	& =-2\det \left( {{D}_{n-1}} \right)-4\det \left( {{D}_{n-2}} \right)+4\left| \begin{matrix}
		{{D}_{n-4}} & \alpha  & \alpha +e  \\
		{{\alpha }^{T}} & 0 & 1  \\
		{{\alpha }^{T}}+2{{e}^{T}} & 2 & 1  \\
	\end{matrix} \right|. 
\end{aligned}
\end{equation*}
Using the property of determinants again, we construct $\det \left( {{D}_{n-2}} \right)$ and do elementary orpertion $R_2\rightarrow R_3$ to simplify the last term into 
\begin{equation*}
      \begin{aligned}
&\left| \begin{matrix}
		{{D}_{n-4}} & \alpha  & \alpha +e  \\
		{{\alpha }^{T}} & 0 & 1  \\
		{{\alpha }^{T}}+2{{e}^{T}} & 2 & 1  \\
	\end{matrix} \right|
=\left| \begin{matrix}
			{{D}_{n-4}} & \alpha  & \alpha +e  \\
			{{\alpha }^{T}} & 0 & 1  \\
			{{\alpha }^{T}}+{{e}^{T}} & 1 & 0  \\
		\end{matrix} \right|+\left| \begin{matrix}
			{{D}_{n-4}} & \alpha  & \alpha +e  \\
			{{\alpha }^{T}} & 0 & 1  \\
			{{e}^{T}} & 1 & 1  \\
		\end{matrix} \right| \\ 
		& =\det \left( {{D}_{n-2}} \right)+\left| \begin{matrix}
			{{D}_{n-4}} & \alpha  & \alpha +e  \\
			{{\alpha }^{T}} & 0 & 1  \\
			{{\alpha }^{T}}+{{e}^{T}} & 1 & 2  \\
		\end{matrix} \right|  \\ &=\det \left( {{D}_{n-2}} \right)+\left| \begin{matrix}
				{{D}_{n-4}} & \alpha  & \alpha +e  \\
				{{\alpha }^{T}} & 0 & 1  \\
				{{\alpha }^{T}}+{{e}^{T}} & 1 & 0  \\
			\end{matrix} \right|+\left| \begin{matrix}
				{{D}_{n-4}} & \alpha  & \alpha +e  \\
				{{\alpha }^{T}} & 0 & 1  \\
				0 & 0 & 2  \\
			\end{matrix} \right| \\ 
			& =2\det \left( {{D}_{n-2}} \right)+2\det \left( {{D}_{n-3}} \right).
	\end{aligned}
\end{equation*}
Hence we have 
\begin{equation*}
	\begin{aligned}
			 \det ({{D}_{n}}) & =-2\det \left( {{D}_{n-1}} \right)-4\det \left( {{D}_{n-2}} \right)+4[2\det \left( {{D}_{n-2}} \right)+2\det \left( {{D}_{n-3}} \right)]\\
& =-2\det \left( {{D}_{n-1}} \right)+4\det \left( {{D}_{n-2}} \right)+8\det \left( {{D}_{n-3}} \right), 
	\end{aligned}
\end{equation*}
i.e., 
\begin{equation}\label{eqn-rr-S4-P4n-3} 
 \det \left( {{D}_{n}} \right)+2\det \left( {{D}_{n-1}} \right)-4\det \left( {{D}_{n-2}} \right)-8\det \left( {{D}_{n-3}} \right)=0.
\end{equation}
Therefore, we have derived the first four-term recurrence relation for the determinant of the distance matrix of trees.


\subsection{The recurrence relation related to $T_n=P_4\circ_{n-2} R$}\label{subsec-S4-P4n-2}

In this subsection, we consider the decomposition $T_n=P_4\circ_{n-2} R$ of $T_n$, and we interchange the labels of $n-2$ and $n-3$, then the partitioned matrix $D_n$ has the following form:

\begin{equation}\label{eqn-S4-P4n-2} 
\det \left( {{D}_{n}} \right)=\left( \begin{matrix}
	{{D}_{n-4}} & \alpha  & \alpha +e & \alpha +e & \alpha +2e  \\
	{{\alpha }^{T}} & 0 & 1 & 1 & 2  \\
	{{\alpha }^{T}}+{{e}^{T}} & 1 & 0 & 2 & 3  \\
	{{\alpha }^{T}}+{{e}^{T}} & 1 & 2 & 0 & 1  \\
	{{\alpha }^{T}}+2{{e}^{T}} & 2 & 3 & 1 & 0  \\
\end{matrix} \right),
\end{equation}
where $\alpha$ is a column vector whose entries are $d_{i,n-3}$ for $i\in [n-4]$. 

We do ${{R}_{2}}-2{{R}_{4}}\rightarrow {{R}_{5}}$ on $\det(D_n)$ and then expand the last row to obtain 
\begin{equation*}
\begin{aligned}
	\det \left( {{D}_{n}} \right)&=\left| \begin{matrix}
		{{D}_{n-4}} & \alpha  & \alpha +e & \alpha +e & \alpha +2e  \\
		{{\alpha }^{T}} & 0 & 1 & 1 & 2  \\
		{{\alpha }^{T}}+{{e}^{T}} & 1 & 0 & 2 & 3  \\
		{{\alpha }^{T}}+{{e}^{T}} & 1 & 2 & 0 & 1  \\
		0 & 0 & 0 & 2 & 0  \\
	\end{matrix} \right|  
=-2\left| \begin{matrix}
		{{D}_{n-4}} & \alpha  & \alpha +e & \alpha +2e  \\
		{{\alpha }^{T}} & 0 & 1 & 2  \\
		{{\alpha }^{T}}+{{e}^{T}} & 1 & 0 & 3  \\
		{{\alpha }^{T}}+{{e}^{T}} & 1 & 2 & 1  \\
	\end{matrix} \right|.  
\end{aligned}
\end{equation*}
We do ${{R}_{3}}-{{R}_{2}}\rightarrow {{R}_{4}}$, then we have 
\begin{equation*}
\begin{aligned}
	\det \left( {{D}_{n}} \right)&=-2\left| \begin{matrix}
		{{D}_{n-4}} & \alpha  & \alpha +e & \alpha +2e  \\
		{{\alpha }^{T}} & 0 & 1 & 2  \\
		{{\alpha }^{T}}+{{e}^{T}} & 1 & 0 & 3  \\
		{{\alpha }^{T}}+2{{e}^{T}} & 2 & 1 & 2  \\
	\end{matrix} \right|
\end{aligned}
\end{equation*}
which coincides with the third line of Equation~(\ref{eqn-S4-P4n-3-4by4}). Hence we get the same recurrence relation as in Equation~(\ref{eqn-rr-S4-P4n-3}).


\subsection{The recurrence relation related to $T_n=K_{1,3}\circ_{n-2} R$}\label{subsec-S4-K13n-2}

In this subsection, we consider the decomposition $T_n=K_{1,3}\circ_{n-2} R$ of $T_n$, and we interchange the labels of $n-2$ and $n-3$, then the partitioned matrix $D_n$ has the following form:

\begin{equation}\label{eqn-S4-K13n-2} 
\det \left( {{D}_{n}} \right)=\left( \begin{matrix}
	{{D}_{n-4}} & \alpha  & \alpha +e & \alpha +2e & \alpha +2e  \\
	{{\alpha }^{T}} & 0 & 1 & 2 & 2  \\
	{{\alpha }^{T}}+{{e}^{T}} & 1 & 0 & 1 & 1  \\
	{{\alpha }^{T}}+2{{e}^{T}} & 2 & 1 & 0 & 2  \\
	{{\alpha }^{T}}+2{{e}^{T}} & 2 & 1 & 2 & 0  \\
\end{matrix} \right), 
\end{equation}
where $\alpha$ is a column vector whose entries are $d_{i,n-3}$ for $i\in [n-4]$. 

We do ${{R}_{2}}-2{{R}_{3}}\rightarrow {{R}_{5}}$ on $\det(D_n)$ and expand it by the last row to get 
\begin{equation*}
\begin{aligned}
	\det \left( {{D}_{n}} \right)& =\left| \begin{matrix}
		{{D}_{n-4}} & \alpha  & \alpha +e & \alpha +2e & \alpha +2e  \\
		{{\alpha }^{T}} & 0 & 1 & 2 & 2  \\
		{{\alpha }^{T}}+{{e}^{T}} & 1 & 0 & 1 & 1  \\
		{{\alpha }^{T}}+2{{e}^{T}} & 2 & 1 & 0 & 2  \\
		0 & 0 & 2 & 2 & 0  \\
	\end{matrix} \right| \\ 
	& =2\left| \begin{matrix}
		{{D}_{n-4}} & \alpha  & \alpha +2e & \alpha +2e  \\
		{{\alpha }^{T}} & 0 & 2 & 2  \\
		{{\alpha }^{T}}+{{e}^{T}} & 1 & 1 & 1  \\
		{{\alpha }^{T}}+2{{e}^{T}} & 2 & 0 & 2  \\
	\end{matrix} \right|-2\left| \begin{matrix}
		{{D}_{n-4}} & \alpha  & \alpha +e & \alpha +2e  \\
		{{\alpha }^{T}} & 0 & 1 & 2  \\
		{{\alpha }^{T}}+{{e}^{T}} & 1 & 0 & 1  \\
		{{\alpha }^{T}}+2{{e}^{T}} & 2 & 1 & 2  \\
	\end{matrix} \right|  
\end{aligned}
\end{equation*}
We do $\frac{1}{2}\left({{C}_{2}}- {{C}_{4}} \right)\rightarrow {{C}_{3}}$ on the first determinant, therefore 
\begin{equation*}
	\begin{aligned}
		\det \left( {{D}_{n}} \right)=2\left| \begin{matrix}
			{{D}_{n-4}} & \alpha  & \alpha +e & \alpha +2e  \\
			{{\alpha }^{T}} & 0 & 1 & 2  \\
			{{\alpha }^{T}}+{{e}^{T}} & 1 & 1 & 1  \\
			{{\alpha }^{T}}+2{{e}^{T}} & 2 & 0 & 2  \\
		\end{matrix} \right|-2\left| \begin{matrix}
			{{D}_{n-4}} & \alpha  & \alpha +e & \alpha +2e  \\
			{{\alpha }^{T}} & 0 & 1 & 2  \\
			{{\alpha }^{T}}+{{e}^{T}} & 1 & 0 & 1  \\
			{{\alpha }^{T}}+2{{e}^{T}} & 2 & 1 & 2  \\
		\end{matrix} \right|  
	\end{aligned}
\end{equation*}
By repeatedly applying the property of determinants to the preceding and following determinants to construct $\det \left( {{D}_{n-1}} \right)$, we obtain 
\begin{equation*}
	\begin{aligned}
\hspace{-3em} \left| \begin{matrix}
			{{D}_{n-4}} & \alpha  & \alpha +e & \alpha +2e  \\
			{{\alpha }^{T}} & 0 & 1 & 2  \\
			{{\alpha }^{T}}+{{e}^{T}} & 1 & 1 & 1  \\
			{{\alpha }^{T}}+2{{e}^{T}} & 2 & 0 & 2  \\
		\end{matrix} \right| 
&= \left| \begin{matrix}
			{{D}_{n-4}} & \alpha  & \alpha +e & \alpha +2e  \\
			{{\alpha }^{T}} & 0 & 1 & 2  \\
			{{\alpha }^{T}}+{{e}^{T}} & 1 & 1 & 1  \\
			{{\alpha }^{T}}+2{{e}^{T}} & 2 & 1 & 0  \\
		\end{matrix} \right|+\left| \begin{matrix}
			{{D}_{n-4}} & \alpha  & \alpha +e & \alpha +2e  \\
			{{\alpha }^{T}} & 0 & 1 & 2  \\
			{{\alpha }^{T}}+{{e}^{T}} & 1 & 1 & 1  \\
			0 & 0 & -1 & 2  \\
		\end{matrix} \right|\\ 
&=  \det(D_{n-1})+\left| \begin{matrix}
		{{D}_{n-4}} & \alpha  & \alpha +2e  \\
		{{\alpha }^{T}} & 0 & 2  \\
		{{\alpha }^{T}}+2{{e}^{T}} & 2 & 0  \\
	\end{matrix} \right|+ \\ &~~~~~~~~~~~~~~~~ 2(\det(D_{n-2})+ \det(D_{n-3}))+\left| \begin{matrix}
			{{D}_{n-4}} & \alpha  & \alpha +2e  \\
			{{\alpha }^{T}} & 0 & 2  \\
			{{\alpha }^{T}}+{{e}^{T}} & 1 & 1  \\
		\end{matrix} \right|,  
	\end{aligned}
\end{equation*}
\begin{equation*}
	\begin{aligned}
		\left| \begin{matrix}
			{{D}_{n-4}} & \alpha  & \alpha +e & \alpha +2e  \\
			{{\alpha }^{T}} & 0 & 1 & 2  \\
			{{\alpha }^{T}}+{{e}^{T}} & 1 & 0 & 1  \\
			{{\alpha }^{T}}+2{{e}^{T}} & 2 & 1 & 2  \\
		\end{matrix} \right|=\det(D_{n-1})+2\det(D_{n-2}). 
	\end{aligned}
\end{equation*}
By doing $C_2\rightarrow C_3$ and $R_2\rightarrow R_3$ to the first determinant, and $C_2\rightarrow C_3$ to the second determinant, we have 
\begin{equation*}
	\begin{aligned}
\left| \begin{matrix}
		{{D}_{n-4}} & \alpha  & \alpha +2e  \\
		{{\alpha }^{T}} & 0 & 2  \\
		{{\alpha }^{T}}+2{{e}^{T}} & 2 & 0  \\
	\end{matrix} \right|
&= 4\left| \begin{matrix}
		{{D}_{n-4}} & \alpha  & \alpha +e  \\
		{{\alpha }^{T}} & 0 & 1  \\
		{{\alpha }^{T}}+{{e}^{T}} & 1 & 1  \\
	\end{matrix} \right|=4(\det(D_{n-2})+\det(D_{n-3})), 
	\end{aligned}
\end{equation*}
\begin{equation*}
	\begin{aligned}
\left| \begin{matrix}
			{{D}_{n-4}} & \alpha  & \alpha +2e  \\
			{{\alpha }^{T}} & 0 & 2  \\
			{{\alpha }^{T}}+{{e}^{T}} & 1 & 1  \\
		\end{matrix} \right|
&= 2\left| \begin{matrix}
		{{D}_{n-4}} & \alpha  & \alpha +e  \\
		{{\alpha }^{T}} & 0 & 1  \\
		{{\alpha }^{T}}+{{e}^{T}} & 1 & 1  \\
	\end{matrix} \right|=2(\det(D_{n-2})+\det(D_{n-3})).
	\end{aligned}
\end{equation*}
Hence we have 
\begin{equation*}
\begin{aligned}
 \det \left( {{D}_{n}} \right) 
& =2[\det(D_{n-1})+8(\det(D_{n-2})+\det(D_{n-3}))] -2[\det(D_{n-1})+2\det(D_{n-2})]\\
&= 12\det \left( {{D}_{n-2}} \right)+16\det \left( {{D}_{n-3}} \right), 
\end{aligned}
\end{equation*}
i.e., 
\begin{equation}\label{eqn-rr-S4-K13n-2}
		\det \left( {{D}_{n}} \right)-12\det \left( {{D}_{n-2}} \right)-16\det \left( {{D}_{n-3}} \right)=0.
\end{equation}
It can be found that Equation~(\ref{eqn-rr-S4-K13n-2}) is a new recurrence relation for the determinant of the distance matrix of trees.


\subsection{The recurrence relation related to $T_n=K_{1,3}\circ_{n-3} R$}\label{subsec-S4-K13n-3}

In this subsection, we consider the decomposition $T_n=K_{1,3}\circ_{n-3} R$ of $T_n$, then the partitioned matrix $D_n$ has the following form:

\begin{equation}\label{eqn-S4-K13n-2}
 {{D}_{n}} =\left( \begin{matrix}
	{{D}_{n-4}} & \alpha  & \alpha +e & \alpha +e & \alpha +e  \\
	{{\alpha }^{T}} & 0 & 1 & 1 & 1  \\
	{{\alpha }^{T}}+{{e}^{T}} & 1 & 0 & 2 & 2  \\
	{{\alpha }^{T}}+{{e}^{T}} & 1 & 2 & 0 & 2  \\
	{{\alpha }^{T}}+{{e}^{T}} & 1 & 2 & 2 & 0  \\
\end{matrix} \right), 
\end{equation}
where $\alpha$ is a column vector whose entries are $d_{i,n-3}$ for $i\in [n-4]$. 

We take $-{{R}_{4}}\rightarrow {{R}_{5}}$, $-{{C}_{4}}\rightarrow {{C}_{5}}$ and expand it by the last row, then  
\begin{equation*}
	\begin{aligned}
\hspace{-3em}	  \det \left( {{D}_{n}} \right)&=\left| \begin{matrix}
				{{D}_{n-4}} & \alpha  & \alpha +e & \alpha +e & 0  \\
				{{\alpha }^{T}} & 0 & 1 & 1 & 0  \\
				{{\alpha }^{T}}+{{e}^{T}} & 1 & 0 & 2 & 0  \\
				{{\alpha }^{T}}+{{e}^{T}} & 1 & 2 & 0 & 2  \\
				0 & 0 & 0 & 2 & -4  \\
			\end{matrix} \right| 
 =-4\det \left( {{D}_{n-1}} \right)-2\left| \begin{matrix}
				{{D}_{n-4}} & \alpha  & \alpha +e & 0  \\
				{{\alpha }^{T}} & 0 & 1 & 0  \\
				{{\alpha }^{T}}+{{e}^{T}} & 1 & 0 & 0  \\
				{{\alpha }^{T}}+{{e}^{T}} & 1 & 2 & 2  \\
			\end{matrix} \right| \\ 
			& =-4\det \left( {{D}_{n-1}} \right)-4\det \left( {{D}_{n-2}} \right), 	
	\end{aligned}
\end{equation*}
i.e., 
\begin{equation}\label{eqn-rr-S4-K13n-3}
		\det \left( {{D}_{n}} \right)+4\det \left( {{D}_{n-1}} \right)+4\det \left( {{D}_{n-2}} \right)=0.
\end{equation}
This three-term recurrence relation has already been found by Du and Yeh~\cite{DuYeh2020}.


\subsection{Summary of recurrence relations related to $S_4$}\label{subsec-S4-summary}

In this section, we considered all the four cases of decompositions $T_n=S_4\circ R$, and we get all the recurrence relations listed in Table~\ref{tab-S4}.

\begin{table}[ht]
	\caption{Recurrence relations related to $T_n=S_4\circ R$}
	\resizebox{\columnwidth}{!}{
	\label{tab-S4}
	\begin{tabular}{|c|c|c|}
		\hline
		decomposition  & recurrence relation & remark \\ \hline
		\begin{tabular}{c}$T_n=P_4\circ_{n-3} R$, or\\ $T_n=P_4\circ_{n-2} R$~~~~\end{tabular}	&  $\det \left( {{D}_{n}} \right)+2\det \left( {{D}_{n-1}} \right)-4\det \left( {{D}_{n-2}} \right)-8\det \left( {{D}_{n-3}} \right)=0$ & four-term\\ \hline
		$T_n=K_{1,3}\circ_{n-3} R$	&  $\det \left( {{D}_{n}} \right)-12\det \left( {{D}_{n-2}} \right)-16\det \left( {{D}_{n-3}} \right)=0$ & four-term  \\ \hline
		$T_n=K_{1,3}\circ_{n-2} R$	&   $\det \left( {{D}_{n}} \right)+4\det \left( {{D}_{n-1}} \right)+4\det \left( {{D}_{n-2}} \right)=0$ & three-term \\ \hline
	\end{tabular}}
\end{table}



\section{Recurrence relations related to $T_n=S_5\circ R$}\label{sec-S5}

In this section, we consider the decomposition $T_n=S_5\circ R$ and recurrence relations related. We may suppose the vertex set of the subtree $S_5$ on five vertices is \begin{center}$[n]\setminus [n-5]=\{n-4,n-3,n-2,n-1,n\}$.\end{center} We know the subtree $S_5$ is a path $P_5$, a star-like tree $T_{5,2}=K(1,1,2)$ or a star $K_{1,4}$, where 
the path $P_5$ has edges \begin{center}$\{n-4,n-3\}$, $\{n-3,n-2\}$, $\{n-2,n-1\}$ and $\{n-1,n\}$,\end{center} 
$T_{5,2}$ has edges \begin{center}$\{n-4,n-2\}$, $\{n-3,n-2\}$, $\{n-2,n-1\}$ and $\{n-1,n\}$,\end{center} and 
$K_{1,4}$ has edges \begin{center}$\{n-4,n-3\}$, $\{n-4,n-2\}$, $\{n-4,n-1\}$ and $\{n-4,n\}$.\end{center}

The decompositions $T_n=P_5\circ_{n} R$ and $T_n=P_5\circ_{n-4} R$ are equivalent, decompositions $T_n=P_5\circ_{n-1} R$ and $T_n=P_5\circ_{n-3} R$ are equivalent, so there are three different decompositions with form $T_n=P_5\circ R$. 

The decompositions $T_n=T_{5,2}\circ_{n-4} R$ and $T_n=T_{5,2}\circ_{n-3} R$ are equivalent, so there are four different decompositions with form $T_n=T_{5,2}\circ R$. 

The decompositions $T_n=K_{1,4}\circ_{n} R$, $T_n=K_{1,4}\circ_{n-1} R$, $T_n=K_{1,4}\circ_{n-2} R$ and $T_n=K_{1,4}\circ_{n-3} R$ are equivalent, so there are two different decompositions with form $T_n=K_{1,4}\circ R$. 

Hence there are nine different decompositions of the form $T_n=S_5\circ R$: 
\begin{center}
  $T_n=P_5\circ_{n-4} R$, $T_n=P_5\circ_{n-3} R$, $T_n=P_5\circ_{n-2} R$,\\ 
  $T_n=T_{5,2}\circ_{n-4} R$, $T_n=T_{5,2}\circ_{n-2} R$, $T_n=T_{5,2}\circ_{n-1} R$, $T_n=T_{5,2}\circ_{n} R$,\\ 
  $T_n=K_{1,4}\circ_{n-4} R$, $T_n=K_{1,4}\circ_{n-3} R$.
\end{center}

For the convenience of calculations in the following subsections, we usually give new labels of vertices of $S_5$ such that label $n-4$ as the common vertex of the decomposition $T_n=S_5\circ R$, and label other vertices of $S_5$ by distances from $n-4$.


\subsection{The recurrence relation related to $T_n=P_5\circ_{n-4} R$}\label{subsec-S5-P5n-4}

In this subsection, we consider the decomposition $T_n=P_5\circ_{n-4} R$ of $T_n$, then the partitioned matrix $D_n$ has the following form:

\begin{equation}\label{eqn-S5-P5n-4}
 {{D}_{n}} =\left( \begin{matrix}
		{{D}_{n-5}} & \alpha  & \alpha +e & \alpha +2e & \alpha +3e & \alpha +4e  \\
		{{\alpha }^{T}} & 0 & 1 & 2 & 3 & 4  \\
		{{\alpha }^{T}}+{{e}^{T}} & 1 & 0 & 1 & 2 & 3  \\
		{{\alpha }^{T}}+2{{e}^{T}} & 2 & 1 & 0 & 1 & 2  \\
		{{\alpha }^{T}}+3{{e}^{T}} & 3 & 2 & 1 & 0 & 1  \\
		{{\alpha }^{T}}+4{{e}^{T}} & 4 & 3 & 2 & 1 & 0  \\
	\end{matrix} \right), 
\end{equation}
where $\alpha$ is a column vector whose entries are $d_{i,n-4}$ for $i\in [n-5]$. 

We take ${{R}_{4}}-2{{R}_{5}}\rightarrow {{R}_{6}}$ and expand it by the last row, then 
\begin{equation*}
	\begin{aligned}
\hspace{-5em}	 \det \left( {{D}_{n}} \right)&=\left| \begin{matrix}
		{{D}_{n-5}} & \alpha  & \alpha +e & \alpha +2e & \alpha +3e & \alpha +4e  \\
		{{\alpha }^{T}} & 0 & 1 & 2 & 3 & 4  \\
		{{\alpha }^{T}}+{{e}^{T}} & 1 & 0 & 1 & 2 & 3  \\
		{{\alpha }^{T}}+2{{e}^{T}} & 2 & 1 & 0 & 1 & 2  \\
		{{\alpha }^{T}}+3{{e}^{T}} & 3 & 2 & 1 & 0 & 1  \\
		0 & 0 & 0 & 0 & 2 & 0  \\
	\end{matrix} \right| 
 =-2\left| \begin{matrix}
		{{D}_{n-5}} & \alpha  & \alpha +e & \alpha +2e & \alpha +4e  \\
		{{\alpha }^{T}} & 0 & 1 & 2 & 4  \\
		{{\alpha }^{T}}+{{e}^{T}} & 1 & 0 & 1 & 3  \\
		{{\alpha }^{T}}+2{{e}^{T}} & 2 & 1 & 0 & 2  \\
		{{\alpha }^{T}}+3{{e}^{T}} & 3 & 2 & 1 & 1  \\
	\end{matrix} \right|. 
	\end{aligned}
\end{equation*}
We now take ${{C}_{3}}-{{C}_{4}}\rightarrow {{C}_{5}}$, and get 
\begin{equation*}
\begin{aligned}
&   \det \left( {{D}_{n}} \right)  =-2\left| \begin{matrix}
	{{D}_{n-5}} & \alpha  & \alpha +e & \alpha +2e & \alpha +3e  \\
	{{\alpha }^{T}} & 0 & 1 & 2 & 3  \\
	{{\alpha }^{T}}+{{e}^{T}} & 1 & 0 & 1 & 2  \\
	{{\alpha }^{T}}+2{{e}^{T}} & 2 & 1 & 0 & 3  \\
	{{\alpha }^{T}}+3{{e}^{T}} & 3 & 2 & 1 & 2  \\
\end{matrix} \right|\\ 
%
 & =-2\left| \begin{matrix}
		{{D}_{n-5}} & \alpha  & \alpha +e & \alpha +2e & \alpha +3e  \\
		{{\alpha }^{T}} & 0 & 1 & 2 & 3  \\
		{{\alpha }^{T}}+{{e}^{T}} & 1 & 0 & 1 & 2  \\
		{{\alpha }^{T}}+2{{e}^{T}} & 2 & 1 & 0 & 1  \\
		{{\alpha }^{T}}+3{{e}^{T}} & 3 & 2 & 1 & 0  \\
	\end{matrix} \right|-2\left| \begin{matrix}
		{{D}_{n-5}} & \alpha  & \alpha +e & \alpha +2e & 0  \\
		{{\alpha }^{T}} & 0 & 1 & 2 & 0  \\
		{{\alpha }^{T}}+{{e}^{T}} & 1 & 0 & 1 & 0  \\
		{{\alpha }^{T}}+2{{e}^{T}} & 2 & 1 & 0 & 2  \\
		{{\alpha }^{T}}+3{{e}^{T}} & 3 & 2 & 1 & 2  \\
	\end{matrix} \right| \\ 
	& =-2\det ({{D}_{n-1}})-4\det ({{D}_{n-2}})+4\left| \begin{matrix}
		{{D}_{n-5}} & \alpha  & \alpha +e & \alpha +2e  \\
		{{\alpha }^{T}} & 0 & 1 & 2  \\
		{{\alpha }^{T}}+{{e}^{T}} & 1 & 0 & 1  \\
		{{\alpha }^{T}}+3{{e}^{T}} & 3 & 2 & 1  \\
	\end{matrix} \right|.  
\end{aligned}
\end{equation*}
Applying the property of determinants of the last determinant, and then taking ${{R}_{2}}-{{R}_{3}}\rightarrow {{R}_{4}}$ on the second determinant and expanding the resulting determiant along the last row, lastly taking $C_2\rightarrow C_3$ on the determinant of order three, we have 
\begin{equation*}
\begin{aligned}
& \hspace{-3em}   \left| \begin{matrix}
        {{D}_{n-5}} & \alpha  & \alpha +e & \alpha +2e  \\
		{{\alpha }^{T}} & 0 & 1 & 2  \\
		{{\alpha }^{T}}+{{e}^{T}} & 1 & 0 & 1  \\
		{{\alpha }^{T}}+3{{e}^{T}} & 3 & 2 & 1  \\
	\end{matrix} \right| 
= \left| \begin{matrix}
		{{D}_{n-5}} & \alpha  & \alpha +e & \alpha +2e  \\
		{{\alpha }^{T}} & 0 & 1 & 2  \\
		{{\alpha }^{T}}+{{e}^{T}} & 1 & 0 & 1  \\
		{{\alpha }^{T}}+2{{e}^{T}} & 2 & 1 & 0  \\
	\end{matrix} \right|   +   \left| \begin{matrix}
		{{D}_{n-5}} & \alpha  & \alpha +e & \alpha +2e  \\
		{{\alpha }^{T}} & 0 & 1 & 2  \\
		{{\alpha }^{T}}+{{e}^{T}} & 1 & 0 & 1  \\
		{{e}^{T}} & 1 & 1 & 1  \\
	\end{matrix} \right|\\ 
& \hspace{-3em}  =   \det ({{D}_{n-2}})+\left| \begin{matrix}
	{{D}_{n-5}} & \alpha  & \alpha +e & \alpha +2e  \\
	{{\alpha }^{T}} & 0 & 1 & 2  \\
	{{\alpha }^{T}}+{{e}^{T}} & 1 & 0 & 1  \\
	0 & 0 & 2 & 2
\end{matrix} \right|
= \det ({{D}_{n-2}})+2\det ({{D}_{n-3}})-2\left| \begin{matrix}
		{{D}_{n-5}} & \alpha  & \alpha +2e  \\
		{{\alpha }^{T}} & 0 & 2  \\
		{{\alpha }^{T}}+{{e}^{T}} & 1 & 1  \\
	\end{matrix} \right|\\
&\hspace{-3em}= \det ({{D}_{n-2}})+2\det ({{D}_{n-3}})-4\left| \begin{matrix}
		{{D}_{n-5}} & \alpha  & \alpha +e  \\
		{{\alpha }^{T}} & 0 & 1  \\
		{{\alpha }^{T}}+{{e}^{T}} & 1 & 1  \\
	\end{matrix} \right|\\
&\hspace{-3em}= \det ({{D}_{n-2}})+2\det ({{D}_{n-3}})-4\left| \begin{matrix}
		{{D}_{n-5}} & \alpha  & \alpha +e  \\
		{{\alpha }^{T}} & 0 & 1  \\
		{{\alpha }^{T}}+{{e}^{T}} & 1 & 0  \\
	\end{matrix} \right|-4\left| \begin{matrix}
		{{D}_{n-5}} & \alpha  & 0  \\
		{{\alpha }^{T}} & 0 & 0  \\
		{{\alpha }^{T}}+{{e}^{T}} & 1 & 1  \\
	\end{matrix} \right|\\
&\hspace{-3em}= \det ({{D}_{n-2}})+2\det ({{D}_{n-3}})-4\det ({{D}_{n-3}})-4\det ({{D}_{n-4}})\\
&\hspace{-3em}= \det ({{D}_{n-2}})-2\det ({{D}_{n-3}})-4\det ({{D}_{n-4}}).
\end{aligned}
\end{equation*}
Hence 
\begin{equation*}
	\begin{aligned}
\det \left( {{D}_{n}} \right)&=-2\det ({{D}_{n-1}})-4\det ({{D}_{n-2}})+4[\det ({{D}_{n-2}})-2\det ({{D}_{n-3}})-4\det ({{D}_{n-4}})]\\ 
			& = -2\det \left( {{D}_{n-1}} \right)-8\det \left( {{D}_{n-3}} \right)-16\det \left( {{D}_{n-4}} \right), 	
	\end{aligned}
\end{equation*}
i.e., 
\begin{equation}\label{eqn-rr-S5-P5n-4}
		\det \left( {{D}_{n}} \right)+2\det \left( {{D}_{n-1}} \right)+8\det \left( {{D}_{n-3}} \right)+16\det \left( {{D}_{n-4}} \right)=0.
\end{equation}

Thus, we have derived the first five-term recurrence relation for the determinant of the distance matrix of trees.


\subsection{The recurrence relation related to $T_n=P_5\circ_{n-3} R$}\label{subsec-S5-P5n-3}

In this subsection, we consider the decomposition $T_n=P_5\circ_{n-3} R$ of $T_n$, and we interchange the labels of $n-3$ and $n-4$, then the partitioned matrix $D_n$ has the following form:

\begin{equation}\label{eqn-S5-P5n-3}
 {{D}_{n}} =\left( \begin{matrix}
    	{{D}_{n-5}} & \alpha  & \alpha +e & \alpha +e & \alpha +2e & \alpha +3e  \\
    	{{\alpha }^{T}} & 0 & 1 & 1 & 2 & 3  \\
    	{{\alpha }^{T}}+{{e}^{T}} & 1 & 0 & 2 & 3 & 4  \\
    	{{\alpha }^{T}}+{{e}^{T}} & 1 & 2 & 0 & 1 & 2  \\
    	{{\alpha }^{T}}+2{{e}^{T}} & 2 & 3 & 1 & 0 & 1  \\
    	{{\alpha }^{T}}+3{{e}^{T}} & 3 & 4 & 2 & 1 & 0  \\
    \end{matrix} \right),
\end{equation}
where $\alpha$ is a column vector whose entries are $d_{i,n-4}$ for $i\in [n-5]$. 

We take ${{R}_{4}}-2{{R}_{5}}\rightarrow {{R}_{6}}$ and expand it by the last row, take ${{C}_{2}}-{{C}_{3}}\rightarrow {{C}_{5}}$ on the first determinant of the second line, then we get 
\begin{equation*}
\begin{aligned}
	& \det \left( {{D}_{n}} \right)=\left| \begin{matrix}
		{{D}_{n-5}} & \alpha  & \alpha +e & \alpha +e & \alpha +2e & \alpha +3e  \\
		{{\alpha }^{T}} & 0 & 1 & 1 & 2 & 3  \\
		{{\alpha }^{T}}+{{e}^{T}} & 1 & 0 & 2 & 3 & 4  \\
		{{\alpha }^{T}}+{{e}^{T}} & 1 & 2 & 0 & 1 & 2  \\
		{{\alpha }^{T}}+2{{e}^{T}} & 2 & 3 & 1 & 0 & 1  \\
		0 & 0 & 0 & 0 & 2 & 0  \\
	\end{matrix} \right| \\ 
	& =-2\left| \begin{matrix}
		{{D}_{n-5}} & \alpha  & \alpha +e & \alpha +e & \alpha +3e  \\
		{{\alpha }^{T}} & 0 & 1 & 1 & 3  \\
		{{\alpha }^{T}}+{{e}^{T}} & 1 & 0 & 2 & 4  \\
		{{\alpha }^{T}}+{{e}^{T}} & 1 & 2 & 0 & 2  \\
		{{\alpha }^{T}}+2{{e}^{T}} & 2 & 3 & 1 & 1  \\
	\end{matrix} \right| =-2\left| \begin{matrix}
		{{D}_{n-5}} & \alpha  & \alpha +e & \alpha +e & \alpha +2e  \\
		{{\alpha }^{T}} & 0 & 1 & 1 & 2  \\
		{{\alpha }^{T}}+{{e}^{T}} & 1 & 0 & 2 & 5  \\
		{{\alpha }^{T}}+{{e}^{T}} & 1 & 2 & 0 & 1  \\
		{{\alpha }^{T}}+2{{e}^{T}} & 2 & 3 & 1 & 0  \\
	\end{matrix} \right|.
\end{aligned}
\end{equation*}
Applying the property of determinants on the last determinant, take $R_2-2R_3\rightarrow R_4$ on the determinant of the second line, then we have 
\begin{equation*}
\begin{aligned}
\hspace{-5em}	&\hspace{-5em} \det \left( {{D}_{n}} \right)=-2\left| \begin{matrix}
		{{D}_{n-5}} & \alpha  & \alpha +e & \alpha +e & \alpha +2e  \\
		{{\alpha }^{T}} & 0 & 1 & 1 & 2  \\
		{{\alpha }^{T}}+{{e}^{T}} & 1 & 0 & 2 & 3  \\
		{{\alpha }^{T}}+{{e}^{T}} & 1 & 2 & 0 & 1  \\
		{{\alpha }^{T}}+2{{e}^{T}} & 2 & 3 & 1 & 0  \\
	\end{matrix} \right|-2\left| \begin{matrix}
		{{D}_{n-5}} & \alpha  & \alpha +e & \alpha +e & 0  \\
		{{\alpha }^{T}} & 0 & 1 & 1 & 0  \\
		{{\alpha }^{T}}+{{e}^{T}} & 1 & 0 & 2 & 2  \\
		{{\alpha }^{T}}+{{e}^{T}} & 1 & 2 & 0 & 0  \\
		{{\alpha }^{T}}+2{{e}^{T}} & 2 & 3 & 1 & 0  \\
	\end{matrix} \right| \\ 
	&\hspace{-5em} =-2\det ({{D}_{n-1}})-4\left| \begin{matrix}
		{{D}_{n-5}} & \alpha  & \alpha +e & \alpha +e  \\
		{{\alpha }^{T}} & 0 & 1 & 1  \\
		{{\alpha }^{T}}+{{e}^{T}} & 1 & 2 & 0  \\
		{{\alpha }^{T}}+2{{e}^{T}} & 2 & 3 & 1  \\
	\end{matrix} \right|= -2\det ({{D}_{n-1}})-4\left| \begin{matrix}
		{{D}_{n-5}} & \alpha  & \alpha +e & \alpha +e  \\
		{{\alpha }^{T}} & 0 & 1 & 1  \\
		{{\alpha }^{T}}+{{e}^{T}} & 1 & 2 & 0  \\
		0 & 0 & 0 & 2  \\
	\end{matrix} \right|\\ 
	&\hspace{-5em} =-2\det ({{D}_{n-1}})-8\left| \begin{matrix}
		{{D}_{n-5}} & \alpha  & \alpha +e  \\
		{{\alpha }^{T}} & 0 & 1  \\
		{{\alpha }^{T}}+{{e}^{T}} & 1 & 2  \\
	\end{matrix} \right| 
 =-2\det ({{D}_{n-1}})-8\left| \begin{matrix}
		{{D}_{n-5}} & \alpha  & \alpha +e  \\
		{{\alpha }^{T}} & 0 & 1  \\
		{{\alpha }^{T}}+{{e}^{T}} & 1 & 0  \\
	\end{matrix} \right| -8\left| \begin{matrix}
		{{D}_{n-5}} & \alpha  & 0  \\
		{{\alpha }^{T}} & 0 & 0  \\
		{{\alpha }^{T}}+{{e}^{T}} & 1 & 2  \\
	\end{matrix} \right| \\
	&\hspace{-5em} =-2\det ({{D}_{n-1}})-8\det(D_{n-3})-16\det(D_{n-4}),
\end{aligned}
\end{equation*}
i.e.
\begin{equation}\label{eqn-rr-S5-P5n-3}
	\det ({{D}_{n}})+2\det ({{D}_{n-1}})+8\det ({{D}_{n-3}})+16\det ({{D}_{n-4}})=0
\end{equation}
This recurrence relation is the same as Equation~(\ref{eqn-rr-S5-P5n-4}).


\subsection{The recurrence relation related to $T_n=P_5\circ_{n-2} R$}\label{subsec-S5-P5n-2}

In this subsection, we consider the decomposition $T_n=P_5\circ_{n-2} R$ of $T_n$, and we interchange the labels of $n-1$ and $n-4$ and then interchange the labels of $n-2$ and $n-4$, then the partitioned matrix $D_n$ has the following form:

\begin{equation}\label{eqn-S5-P5n-2}
{{D}_{n}}=\left( \begin{matrix}
	{{D}_{n-5}} & \alpha  & \alpha +e & \alpha +e & \alpha +2e & \alpha +2e  \\
	{{\alpha }^{T}} & 0 & 1 & 1 & 2 & 2  \\
	{{\alpha }^{T}}+{{e}^{T}} & 1 & 0 & 2 & 1 & 3  \\
	{{\alpha }^{T}}+{{e}^{T}} & 1 & 2 & 0 & 3 & 1  \\
	{{\alpha }^{T}}+2{{e}^{T}} & 2 & 1 & 3 & 0 & 4  \\
	{{\alpha }^{T}}+2{{e}^{T}} & 2 & 3 & 1 & 4 & 0  \\
\end{matrix} \right), 
\end{equation}
where $\alpha$ is a column vector whose entries are $d_{i,n-4}$ for $i\in [n-5]$.

We take ${{R}_{2}}-2{{R}_{4}}\rightarrow {{R}_{6}}$ and expand it by the last row, swap ${{C}_{5}}$ and ${{C}_{4}}$ of the second determinant of the first line and then take $\frac{1}{2}({{C}_{2}}- {{C}_{5}})\rightarrow {{C}_{4}}$ on the resulting determinant, apply the property of determinants on the second determinant of the second line for its fourth column, expand the second determinant of the third line along the fourth column, take ${{C}_{2}}-{{C}_{3}}\rightarrow {{C}_{4}}$ on the two determinants of the fourth line and take ${{R}_{2}}-{{R}_{3}}\rightarrow {{R}_{4}}$ on the second determinant of the same line, then we get 

\begin{equation*}
	\begin{aligned}
\hspace{-5em}		 \det \left( {{D}_{n}} \right)&=\left| \begin{matrix}
			{{D}_{n-5}} & \alpha  & \alpha +e & \alpha +e & \alpha +2e & \alpha +2e  \\
			{{\alpha }^{T}} & 0 & 1 & 1 & 2 & 2  \\
			{{\alpha }^{T}}+{{e}^{T}} & 1 & 0 & 2 & 1 & 3  \\
			{{\alpha }^{T}}+{{e}^{T}} & 1 & 2 & 0 & 3 & 1  \\
			{{\alpha }^{T}}+2{{e}^{T}} & 2 & 1 & 3 & 0 & 4  \\
			0 & 0 & 0 & 2 & 0 & 0  \\
		\end{matrix} \right| 
 =2\left| \begin{matrix}
			{{D}_{n-5}} & \alpha  & \alpha +e & \alpha +2e & \alpha +2e  \\
			{{\alpha }^{T}} & 0 & 1 & 2 & 2  \\
			{{\alpha }^{T}}+{{e}^{T}} & 1 & 0 & 1 & 3  \\
			{{\alpha }^{T}}+{{e}^{T}} & 1 & 2 & 3 & 1  \\
			{{\alpha }^{T}}+2{{e}^{T}} & 2 & 1 & 0 & 4  \\
		\end{matrix} \right| \\
& =-2\left| \begin{matrix}
       	{{D}_{n-5}} & \alpha  & \alpha +e & \alpha +2e & \alpha +2e  \\
    	{{\alpha }^{T}} & 0 & 1 & 2 & 2  \\
    	{{\alpha }^{T}}+{{e}^{T}} & 1 & 0 & 3 & 1  \\
    	{{\alpha }^{T}}+{{e}^{T}} & 1 & 2 & 1 & 3  \\
    	{{\alpha }^{T}}+2{{e}^{T}} & 2 & 1 & 4 & 0  \\
        \end{matrix} \right|=-2\left| \begin{matrix}
    	{{D}_{n-5}} & \alpha  & \alpha +e & \alpha +e & \alpha +2e  \\
    	{{\alpha }^{T}} & 0 & 1 & 1 & 2  \\
    	{{\alpha }^{T}}+{{e}^{T}} & 1 & 0 & 3 & 1  \\
    	{{\alpha }^{T}}+{{e}^{T}} & 1 & 2 & 0 & 3  \\
    	{{\alpha }^{T}}+2{{e}^{T}} & 2 & 1 & 5 & 0  \\
      \end{matrix} \right| \\
& =-2\left| \begin{matrix}
			{{D}_{n-5}} & \alpha  & \alpha +e & \alpha +e & \alpha +2e  \\
			{{\alpha }^{T}} & 0 & 1 & 1 & 2  \\
			{{\alpha }^{T}}+{{e}^{T}} & 1 & 0 & 2 & 1  \\
			{{\alpha }^{T}}+{{e}^{T}} & 1 & 2 & 0 & 3  \\
			{{\alpha }^{T}}+2{{e}^{T}} & 2 & 1 & 3 & 0  \\
		\end{matrix} \right|-2\left| \begin{matrix}
			{{D}_{n-5}} & \alpha  & \alpha +e & 0 & \alpha +2e  \\
			{{\alpha }^{T}} & 0 & 1 & 0 & 2  \\
			{{\alpha }^{T}}+{{e}^{T}} & 1 & 0 & 1 & 1  \\
			{{\alpha }^{T}}+{{e}^{T}} & 1 & 2 & 0 & 3  \\
			{{\alpha }^{T}}+2{{e}^{T}} & 2 & 1 & 2 & 0  \\
		\end{matrix} \right| \\ 
		& =-2\det ({{D}_{n-1}})+4\left| \begin{matrix}
			{{D}_{n-5}} & \alpha  & \alpha +e & \alpha +2e  \\
			{{\alpha }^{T}} & 0 & 1 & 2  \\
			{{\alpha }^{T}}+{{e}^{T}} & 1 & 0 & 1  \\
			{{\alpha }^{T}}+{{e}^{T}} & 1 & 2 & 3  \\
		\end{matrix} \right|+2\left| \begin{matrix}
			{{D}_{n-5}} & \alpha  & \alpha +e & \alpha +2e  \\
			{{\alpha }^{T}} & 0 & 1 & 2  \\
			{{\alpha }^{T}}+{{e}^{T}} & 1 & 2 & 3  \\
			{{\alpha }^{T}}+2{{e}^{T}} & 2 & 1 & 0  \\
		\end{matrix} \right| \\
& =-2\det ({{D}_{n-1}})+4\left| \begin{matrix}
			{{D}_{n-5}} & \alpha  & \alpha +e & \alpha +e  \\
			{{\alpha }^{T}} & 0 & 1 & 1  \\
			{{\alpha }^{T}}+{{e}^{T}} & 1 & 0 & 2  \\
			{{\alpha }^{T}}+{{e}^{T}} & 1 & 2 & 2  \\
		\end{matrix} \right|+2\left| \begin{matrix}
			{{D}_{n-5}} & \alpha  & \alpha +e & \alpha +e  \\
			{{\alpha }^{T}} & 0 & 1 & 1  \\
			{{\alpha }^{T}}+{{e}^{T}} & 1 & 2 & 2  \\
			{{\alpha }^{T}}+{{e}^{T}} & 1 & 0 & 0  \\
		\end{matrix} \right|. 
	\end{aligned}
\end{equation*}	
We know 
\begin{equation*}
	\begin{aligned}
\left| \begin{matrix}
			{{D}_{n-5}} & \alpha  & \alpha +e & \alpha +e  \\
			{{\alpha }^{T}} & 0 & 1 & 1  \\
			{{\alpha }^{T}}+{{e}^{T}} & 1 & 0 & 2  \\
			{{\alpha }^{T}}+{{e}^{T}} & 1 & 2 & 2  \\
		\end{matrix} \right|=\left| \begin{matrix}
			{{D}_{n-5}} & \alpha  & \alpha +e & \alpha +e +0 \\
			{{\alpha }^{T}} & 0 & 1 & 1 +0 \\
			{{\alpha }^{T}}+{{e}^{T}} & 1 & 0 & 2 +0 \\
			{{\alpha }^{T}}+{{e}^{T}} & 1 & 2 & 0+ 2  \\
		\end{matrix} \right|= \det(D_{n-2})+2\det(D_{n-3}).
	\end{aligned}
\end{equation*}
Applying the property of determinants, we have
\begin{equation*}
	\begin{aligned}
\left| \begin{matrix}
			{{D}_{n-5}} & \alpha  & \alpha +e & \alpha +e  \\
			{{\alpha }^{T}} & 0 & 1 & 1  \\
			{{\alpha }^{T}}+{{e}^{T}} & 1 & 2 & 2  \\
			{{\alpha }^{T}}+{{e}^{T}} & 1 & 0 & 0  \\
		\end{matrix} \right| 
&= -\left| \begin{matrix}
			{{D}_{n-5}} & \alpha  & \alpha +e & \alpha +e  \\
			{{\alpha }^{T}} & 0 & 1 & 1  \\
			{{\alpha }^{T}}+{{e}^{T}} & 1 & 0 & 0  \\
			{{\alpha }^{T}}+{{e}^{T}} & 1 & 2 & 2\\
		\end{matrix} \right|\\
& = -\left| \begin{matrix}
			{{D}_{n-5}} & \alpha  & \alpha +e & \alpha +e  \\
			{{\alpha }^{T}} & 0 & 1 & 1  \\
			{{\alpha }^{T}}+{{e}^{T}} & 1 & 0 & 2  \\
			{{\alpha }^{T}}+{{e}^{T}} & 1 & 2 & 2  \\
		\end{matrix} \right|-\left| \begin{matrix}
		{{D}_{n-5}} & \alpha  & \alpha +e & 0  \\
		{{\alpha }^{T}} & 0 & 1 & 0  \\
		{{\alpha }^{T}}+{{e}^{T}} & 1 & 0 & -2  \\
		{{\alpha }^{T}}+{{e}^{T}} & 1 & 2 & 0  \\
	\end{matrix} \right|,
	\end{aligned}
\end{equation*}
\begin{equation*}
	\begin{aligned}
\left| \begin{matrix}
		{{D}_{n-5}} & \alpha  & \alpha +e & 0  \\
		{{\alpha }^{T}} & 0 & 1 & 0  \\
		{{\alpha }^{T}}+{{e}^{T}} & 1 & 0 & -2  \\
		{{\alpha }^{T}}+{{e}^{T}} & 1 & 2 & 0  \\
	\end{matrix} \right|=2\left| \begin{matrix}
		{{D}_{n-5}} & \alpha  & \alpha +e \\
		{{\alpha }^{T}} & 0 & 1 \\
		{{\alpha }^{T}}+{{e}^{T}} & 1 & 2 \\
	\end{matrix} \right|=2\det(D_{n-3})+4\det(D_{n-4}).
	\end{aligned}
\end{equation*}
Hence we have 
\begin{equation*}
\begin{aligned}
	 \det \left( {{D}_{n}} \right)&=-2\det ({{D}_{n-1}})+2[\det(D_{n-2})+2\det(D_{n-3})]-2[2\det(D_{n-3})+4\det(D_{n-4})]\\ 
	& =-2\det ({{D}_{n-1}})+2\det ({{D}_{n-2}})-8\det ({{D}_{n-4}}), 
\end{aligned}
\end{equation*}
i.e.,
\begin{equation}\label{3-17}
    \det ({{D}_{n}})+2\det ({{D}_{n-1}})-2\det ({{D}_{n-2}})+8\det ({{D}_{n-4}})=0.
\end{equation}	
Therefore, we have derived the second five-term recurrence relation for the determinant of the distance matrix of trees.


\subsection{The recurrence relation related to $T_n=T_{5,2}\circ_{n-4} R$}\label{subsec-S5-T52n-4}

In this subsection, we consider the decomposition $T_n=T_{5,2}\circ_{n-4} R$ of $T_n$, and we interchange the labels of $n-2$ and $n-3$, then the partitioned matrix $D_n$ has the following form:
\begin{equation}\label{eqn-S5-T52n-4}
 {{D}_{n}} =\left( \begin{matrix}
	{{D}_{n-5}} & \alpha  & \alpha +e & \alpha +2e & \alpha +2e & \alpha +3e  \\
	{{\alpha }^{T}} & 0 & 1 & 2 & 2 & 3  \\
	{{\alpha }^{T}}+{{e}^{T}} & 1 & 0 & 1 & 1 & 2  \\
	{{\alpha }^{T}}+2{{e}^{T}} & 2 & 1 & 0 & 2 & 3  \\
	{{\alpha }^{T}}+2{{e}^{T}} & 2 & 1 & 2 & 0 & 1  \\
	{{\alpha }^{T}}+3{{e}^{T}} & 3 & 2 & 3 & 1 & 0  \\
\end{matrix} \right), 
\end{equation}
where $\alpha$ is a column vector whose entries are $d_{i,n-4}$ for $i\in [n-5]$. 

We take ${{R}_{3}}-2{{R}_{5}}\rightarrow {{R}_{6}}$ and expand it by the last row, take ${{C}_{3}}-{{C}_{4}}\rightarrow {{C}_{5}}$ on the second determinant of the first line, then we have 
\begin{equation*}
\begin{aligned}
    &\hspace{-6em}    \det \left( {{D}_{n}} \right)=\left| \begin{matrix}
		{{D}_{n-5}} & \alpha  & \alpha +e & \alpha +2e & \alpha +2e & \alpha +3e  \\
		{{\alpha }^{T}} & 0 & 1 & 2 & 2 & 3  \\
		{{\alpha }^{T}}+{{e}^{T}} & 1 & 0 & 1 & 1 & 2  \\
		{{\alpha }^{T}}+2{{e}^{T}} & 2 & 1 & 0 & 2 & 3  \\
		{{\alpha }^{T}}+2{{e}^{T}} & 2 & 1 & 2 & 0 & 1  \\
		0 & 0 & 0 & 0 & 2 & 0  \\
	\end{matrix} \right| 
 =-2\left| \begin{matrix}
		{{D}_{n-5}} & \alpha  & \alpha +e & \alpha +2e & \alpha +3e  \\
		{{\alpha }^{T}} & 0 & 1 & 2 & 3  \\
		{{\alpha }^{T}}+{{e}^{T}} & 1 & 0 & 1 & 2  \\
		{{\alpha }^{T}}+2{{e}^{T}} & 2 & 1 & 0 & 3  \\
		{{\alpha }^{T}}+2{{e}^{T}} & 2 & 1 & 2 & 1  \\
	\end{matrix} \right| \\
&\hspace{-6em}=-2\left| \begin{matrix}
		{{D}_{n-5}} & \alpha  & \alpha +e & \alpha +2e & \alpha +2e  \\
		{{\alpha }^{T}} & 0 & 1 & 2 & 2  \\
		{{\alpha }^{T}}+{{e}^{T}} & 1 & 0 & 1 & 1  \\
		{{\alpha }^{T}}+2{{e}^{T}} & 2 & 1 & 0 & 4  \\
		{{\alpha }^{T}}+2{{e}^{T}} & 2 & 1 & 2 & 0  \\
	\end{matrix} \right| 
 =-2\det ({{D}_{n-1}})-2\left| \begin{matrix}
		{{D}_{n-5}} & \alpha  & \alpha +e & \alpha +2e & 0  \\
		{{\alpha }^{T}} & 0 & 1 & 2 & 0  \\
		{{\alpha }^{T}}+{{e}^{T}} & 1 & 0 & 1 & 0  \\
		{{\alpha }^{T}}+2{{e}^{T}} & 2 & 1 & 0 & 2  \\
		{{\alpha }^{T}}+2{{e}^{T}} & 2 & 1 & 2 & 0  \\
	\end{matrix} \right| \\ 
	&\hspace{-6em} =-2\det ({{D}_{n-1}})+4\left| \begin{matrix}
		{{D}_{n-5}} & \alpha  & \alpha +e & \alpha +2e  \\
		{{\alpha }^{T}} & 0 & 1 & 2  \\
		{{\alpha }^{T}}+{{e}^{T}} & 1 & 0 & 1  \\
		{{\alpha }^{T}}+2{{e}^{T}} & 2 & 1 & 2  \\
	\end{matrix} \right| 
 =-2\det ({{D}_{n-1}})+4\det ({{D}_{n-2}})+8\det ({{D}_{n-3}}) , 
\end{aligned}
\end{equation*}
i.e., 
\begin{equation}\label{eqn-rr-S5-T52n-4}
\det ({{D}_{n}})+2\det ({{D}_{n-1}})-4\det ({{D}_{n-2}})-8\det ({{D}_{n-3}})=0.
\end{equation}
This recurrence relation is the same as Equation~(\ref{eqn-rr-S4-P4n-3}).


\subsection{The recurrence relation related to $T_n=T_{5,2}\circ_{n-2} R$}\label{subsec-S5-T52n-2}

In this subsection, we consider the decomposition $T_n=T_{5,2}\circ_{n-2} R$ of $T_n$, and we interchange the labels of $n-2$ and $n-4$, then the partitioned matrix $D_n$ has the following form:

\begin{equation}\label{eqn-S5-T52n-2}
 {{D}_{n}} =\left( \begin{matrix}
	{{D}_{n-5}} & \alpha  & \alpha +e & \alpha +e & \alpha +e & \alpha +2e  \\
	{{\alpha }^{T}} & 0 & 1 & 1 & 1 & 2  \\
	{{\alpha }^{T}}+{{e}^{T}} & 1 & 0 & 2 & 2 & 3  \\
	{{\alpha }^{T}}+{{e}^{T}} & 1 & 2 & 0 & 2 & 3  \\
	{{\alpha }^{T}}+{{e}^{T}} & 1 & 2 & 2 & 0 & 1  \\
	{{\alpha }^{T}}+2{{e}^{T}} & 2 & 3 & 3 & 1 & 0  \\
\end{matrix} \right), 
\end{equation}
where $\alpha$ is a column vector whose entries are $d_{i,n-4}$ for $i\in [n-5]$. 

We take ${{R}_{2}}-2{{R}_{5}}\rightarrow {{R}_{6}}$ and expand it by the last row, take ${{C}_{2}}-{{C}_{4}}\rightarrow {{C}_{5}}$ on the second determinant of the first line, then we have 
\begin{equation*}
  \begin{aligned}
\hspace{-4em}	 \det \left( {{D}_{n}} \right)&=\left| \begin{matrix}
		{{D}_{n-5}} & \alpha  & \alpha +e & \alpha +e & \alpha +e & \alpha +2e  \\
		{{\alpha }^{T}} & 0 & 1 & 1 & 1 & 2  \\
		{{\alpha }^{T}}+{{e}^{T}} & 1 & 0 & 2 & 2 & 3  \\
		{{\alpha }^{T}}+{{e}^{T}} & 1 & 2 & 0 & 2 & 3  \\
		{{\alpha }^{T}}+{{e}^{T}} & 1 & 2 & 2 & 0 & 1  \\
		0 & 0 & 0 & 0 & 2 & 0  \\
	\end{matrix} \right| 
 =-2\left| \begin{matrix}
		{{D}_{n-5}} & \alpha  & \alpha +e & \alpha +e & \alpha +2e  \\
		{{\alpha }^{T}} & 0 & 1 & 1 & 2  \\
		{{\alpha }^{T}}+{{e}^{T}} & 1 & 0 & 2 & 3  \\
		{{\alpha }^{T}}+{{e}^{T}} & 1 & 2 & 0 & 3  \\
		{{\alpha }^{T}}+{{e}^{T}} & 1 & 2 & 2 & 1  \\
	\end{matrix} \right|  \\ 
&=-2\left| \begin{matrix}
		{{D}_{n-5}} & \alpha  & \alpha +e & \alpha +e & \alpha +e  \\
		{{\alpha }^{T}} & 0 & 1 & 1 & 1  \\
		{{\alpha }^{T}}+{{e}^{T}} & 1 & 0 & 2 & 2  \\
		{{\alpha }^{T}}+{{e}^{T}} & 1 & 2 & 0 & 4  \\
		{{\alpha }^{T}}+{{e}^{T}} & 1 & 2 & 2 & 0  \\
	\end{matrix} \right| 
 =-2\det ({{D}_{n-1}})-2\left| \begin{matrix}
		{{D}_{n-5}} & \alpha  & \alpha +e & \alpha +e & 0  \\
		{{\alpha }^{T}} & 0 & 1 & 1 & 0  \\
		{{\alpha }^{T}}+{{e}^{T}} & 1 & 0 & 2 & 0  \\
		{{\alpha }^{T}}+{{e}^{T}} & 1 & 2 & 0 & 2  \\
		{{\alpha }^{T}}+{{e}^{T}} & 1 & 2 & 2 & 0  \\
	\end{matrix} \right| \\ 
	& =-2\det ({{D}_{n-1}})+4\left| \begin{matrix}
		{{D}_{n-5}} & \alpha  & \alpha +e & \alpha +e  \\
		{{\alpha }^{T}} & 0 & 1 & 1  \\
		{{\alpha }^{T}}+{{e}^{T}} & 1 & 0 & 2  \\
		{{\alpha }^{T}}+{{e}^{T}} & 1 & 2 & 2  \\
	\end{matrix} \right| 
 =-2\det ({{D}_{n-1}})+4\det ({{D}_{n-2}})+8\det ({{D}_{n-3}}), 
\end{aligned}
\end{equation*}
i.e., 
\begin{equation}\label{eqn-rr-S5-T52n-2}
\det ({{D}_{n}})+2\det ({{D}_{n-1}})-4\det ({{D}_{n-2}})-8\det ({{D}_{n-3}})=0. 
\end{equation}
This recurrence relation is the same as Equation~(\ref{eqn-rr-S4-P4n-3}) and Equation~(\ref{eqn-rr-S5-T52n-4}).


\subsection{The recurrence relation related to $T_n=T_{5,2}\circ_{n-1} R$}\label{subsec-S5-T52n-1}

In this subsection, we consider the decomposition $T_n=T_{5,2}\circ_{n-1} R$ of $T_n$, and we interchange the labels of $n-1$ and $n-4$ and interchange the labels of $n-3$ and $n$, then the partitioned matrix $D_n$ has the following form:
\begin{equation}\label{eqn-S5-T52n-1}
 {{D}_{n}} =\left( \begin{matrix}
	{{D}_{n-5}} & \alpha  & \alpha +e & \alpha +e & \alpha +2e & \alpha +2e  \\
	{{\alpha }^{T}} & 0 & 1 & 1 & 2 & 2  \\
	{{\alpha }^{T}}+{{e}^{T}} & 1 & 0 & 2 & 3 & 3  \\
	{{\alpha }^{T}}+{{e}^{T}} & 1 & 2 & 0 & 1 & 1  \\
	{{\alpha }^{T}}+2{{e}^{T}} & 2 & 3 & 1 & 0 & 2  \\
	{{\alpha }^{T}}+2{{e}^{T}} & 2 & 3 & 1 & 2 & 0  \\
    \end{matrix} \right),
\end{equation}
where $\alpha$ is a column vector whose entries are $d_{i,n-4}$ for $i\in [n-5]$.

We take ${{R}_{2}}-2{{R}_{4}}\rightarrow {{R}_{6}}$ and expand it by the last row, take ${{C}_{2}}-{{C}_{3}}\rightarrow {{C}_{4}}$ on the first determinant of the second line, then we have 
\begin{equation*}
\begin{aligned}
&	 \det \left( {{D}_{n}} \right)=\left| \begin{matrix}
		{{D}_{n-5}} & \alpha  & \alpha +e & \alpha +e & \alpha +2e & \alpha +2e  \\
		{{\alpha }^{T}} & 0 & 1 & 1 & 2 & 2  \\
		{{\alpha }^{T}}+{{e}^{T}} & 1 & 0 & 2 & 3 & 3  \\
		{{\alpha }^{T}}+{{e}^{T}} & 1 & 2 & 0 & 1 & 1  \\
		{{\alpha }^{T}}+2{{e}^{T}} & 2 & 3 & 1 & 0 & 2  \\
		0 & 0 & 0 & 2 & 2 & 0  \\
	\end{matrix} \right| \\ 
	& =2\left| \begin{matrix}
		{{D}_{n-5}} & \alpha  & \alpha +e & \alpha +2e & \alpha +2e  \\
		{{\alpha }^{T}} & 0 & 1 & 2 & 2  \\
		{{\alpha }^{T}}+{{e}^{T}} & 1 & 0 & 3 & 3  \\
		{{\alpha }^{T}}+{{e}^{T}} & 1 & 2 & 1 & 1  \\
		{{\alpha }^{T}}+2{{e}^{T}} & 2 & 3 & 0 & 2  \\
	\end{matrix} \right|-2\left| \begin{matrix}
		{{D}_{n-5}} & \alpha  & \alpha +e & \alpha +e & \alpha +2e  \\
		{{\alpha }^{T}} & 0 & 1 & 1 & 2  \\
		{{\alpha }^{T}}+{{e}^{T}} & 1 & 0 & 2 & 3  \\
		{{\alpha }^{T}}+{{e}^{T}} & 1 & 2 & 0 & 1  \\
		{{\alpha }^{T}}+2{{e}^{T}} & 2 & 3 & 1 & 2  \\
	\end{matrix} \right|  \\
&=2\left| \begin{matrix}
		{{D}_{n-5}} & \alpha  & \alpha +e & \alpha +e & \alpha +2e  \\
		{{\alpha }^{T}} & 0 & 1 & 1 & 2  \\
		{{\alpha }^{T}}+{{e}^{T}} & 1 & 0 & 4 & 3  \\
		{{\alpha }^{T}}+{{e}^{T}} & 1 & 2 & 0 & 1  \\
		{{\alpha }^{T}}+2{{e}^{T}} & 2 & 3 & -1 & 2  \\
	\end{matrix} \right|-2\det ({{D}_{n-1}})-4\det ({{D}_{n-2}}).
\end{aligned}
\end{equation*}
Applying the property of determinants on the last determinant for its fifth column, then we get 
\begin{equation*}
\begin{aligned}
&  \left| \begin{matrix}
		{{D}_{n-5}} & \alpha  & \alpha +e & \alpha +e & \alpha +2e  \\
		{{\alpha }^{T}} & 0 & 1 & 1 & 2  \\
		{{\alpha }^{T}}+{{e}^{T}} & 1 & 0 & 4 & 3  \\
		{{\alpha }^{T}}+{{e}^{T}} & 1 & 2 & 0 & 1  \\
		{{\alpha }^{T}}+2{{e}^{T}} & 2 & 3 & -1 & 2  \\
	\end{matrix} \right| \\
&	=\left| \begin{matrix}
		{{D}_{n-5}} & \alpha  & \alpha +e & \alpha +e & \alpha +2e  \\
		{{\alpha }^{T}} & 0 & 1 & 1 & 2  \\
		{{\alpha }^{T}}+{{e}^{T}} & 1 & 0 & 4 & 3  \\
		{{\alpha }^{T}}+{{e}^{T}} & 1 & 2 & 0 & 1  \\
		{{\alpha }^{T}}+2{{e}^{T}} & 2 & 3 & -1 & 0  \\
	\end{matrix} \right| 
 +\left| \begin{matrix}
		{{D}_{n-5}} & \alpha  & \alpha +e & \alpha +e & 0  \\
		{{\alpha }^{T}} & 0 & 1 & 1 & 0  \\
		{{\alpha }^{T}}+{{e}^{T}} & 1 & 0 & 4 & 0  \\
		{{\alpha }^{T}}+{{e}^{T}} & 1 & 2 & 0 & 0  \\
		{{\alpha }^{T}}+2{{e}^{T}} & 2 & 3 & -1 & 2  \\
	\end{matrix} \right| \\ 
	& = \det ({{D}_{n-1}})+\left| \begin{matrix}
		{{D}_{n-5}} & \alpha  & \alpha +e & 0 & \alpha +2e  \\
		{{\alpha }^{T}} & 0 & 1 & 0 & 2  \\
		{{\alpha }^{T}}+{{e}^{T}} & 1 & 0 & 2 & 3  \\
		{{\alpha }^{T}}+{{e}^{T}} & 1 & 2 & 0 & 1  \\
		{{\alpha }^{T}}+2{{e}^{T}} & 2 & 3 & -2 & 0  \\
	\end{matrix} \right|  +2\left| \begin{matrix}
		{{D}_{n-5}} & \alpha  & \alpha +e & \alpha +e  \\
		{{\alpha }^{T}} & 0 & 1 & 1  \\
		{{\alpha }^{T}}+{{e}^{T}} & 1 & 0 & 4  \\
		{{\alpha }^{T}}+{{e}^{T}} & 1 & 2 & 0  \\ 
	\end{matrix} \right|,
\end{aligned}
\end{equation*}
\begin{equation*}
\begin{aligned}
\hspace{-5em} \left| \begin{matrix}
		{{D}_{n-5}} & \alpha  & \alpha +e & 0 & \alpha +2e  \\
		{{\alpha }^{T}} & 0 & 1 & 0 & 2  \\
		{{\alpha }^{T}}+{{e}^{T}} & 1 & 0 & 2 & 3  \\
		{{\alpha }^{T}}+{{e}^{T}} & 1 & 2 & 0 & 1  \\
		{{\alpha }^{T}}+2{{e}^{T}} & 2 & 3 & -2 & 0  \\
	\end{matrix} \right|=2\left| \begin{matrix}
		{{D}_{n-5}} & \alpha  & \alpha +e  & \alpha +2e  \\
		{{\alpha }^{T}} & 0 & 1  & 2  \\
		{{\alpha }^{T}}+{{e}^{T}} & 1 & 0  & 3  \\
		{{\alpha }^{T}}+{{e}^{T}} & 1 & 2  & 1 
	\end{matrix} \right| -2\left| \begin{matrix}
		{{D}_{n-5}} & \alpha  & \alpha +e  & \alpha +2e  \\
		{{\alpha }^{T}} & 0 & 1  & 2  \\ 
		{{\alpha }^{T}}+{{e}^{T}} & 1 & 2  & 1  \\
		{{\alpha }^{T}}+2{{e}^{T}} & 2 & 3  & 0  \\
	\end{matrix} \right|,
\end{aligned}
\end{equation*}
\begin{equation*}
\begin{aligned}
&\hspace{-4em} \left| \begin{matrix}
		{{D}_{n-5}} & \alpha  & \alpha +e  & \alpha +2e  \\
		{{\alpha }^{T}} & 0 & 1  & 2  \\
		{{\alpha }^{T}}+{{e}^{T}} & 1 & 0  & 3  \\
		{{\alpha }^{T}}+{{e}^{T}} & 1 & 2  & 1 
	\end{matrix} \right|= \left| \begin{matrix}
		{{D}_{n-5}} & \alpha  & \alpha +e  & \alpha +e  \\
		{{\alpha }^{T}} & 0 & 1  & 1  \\
		{{\alpha }^{T}}+{{e}^{T}} & 1 & 0  & 4  \\
		{{\alpha }^{T}}+{{e}^{T}} & 1 & 2  & 0 
	\end{matrix} \right|= \det(D_{n-2})+\left| \begin{matrix}
		{{D}_{n-5}} & \alpha  & \alpha +e & 0  \\
		{{\alpha }^{T}} & 0 & 1 & 0  \\
		{{\alpha }^{T}}+{{e}^{T}} & 1 & 0 & 2  \\
		{{\alpha }^{T}}+{{e}^{T}} & 1 & 2 & 0  \\ 
	\end{matrix} \right|\\
&\hspace{-4em}= \det(D_{n-2})-2\left| \begin{matrix}
		{{D}_{n-5}} & \alpha  & \alpha +e   \\
		{{\alpha }^{T}} & 0 & 1   \\ 
		{{\alpha }^{T}}+{{e}^{T}} & 1 & 2  \\ 
	\end{matrix} \right|= \det(D_{n-2})-2\det(D_{n-3})-4\det(D_{n-4}), 
\end{aligned}
\end{equation*}
\begin{equation*}
\begin{aligned}
&\hspace{-3em} \left| \begin{matrix}
		{{D}_{n-5}} & \alpha  & \alpha +e  & \alpha +2e  \\
		{{\alpha }^{T}} & 0 & 1  & 2  \\ 
		{{\alpha }^{T}}+{{e}^{T}} & 1 & 2  & 1  \\
		{{\alpha }^{T}}+2{{e}^{T}} & 2 & 3  & 0  \\
	\end{matrix} \right|= \left| \begin{matrix}
		{{D}_{n-5}} & \alpha  & \alpha +e  & \alpha +e  \\
		{{\alpha }^{T}} & 0 & 1  & 1  \\ 
		{{\alpha }^{T}}+{{e}^{T}} & 1 & 2  & 0  \\
		{{\alpha }^{T}}+2{{e}^{T}} & 2 & 3  & -1  \\
	\end{matrix} \right| = \left| \begin{matrix}
		{{D}_{n-5}} & \alpha  & \alpha +e  & \alpha +e  \\
		{{\alpha }^{T}} & 0 & 1  & 1  \\ 
		{{\alpha }^{T}}+{{e}^{T}} & 1 & 2  & 0  \\
		{{\alpha }^{T}}+{{e}^{T}} & 1 & 2  & 0  \\
	\end{matrix} \right|=0,  
\end{aligned}
\end{equation*}
Hence we have 
\begin{equation*}
\begin{aligned}
\hspace{-4em} \det \left( {{D}_{n}} \right)&=2[\det(D_{n-1})+4(\det(D_{n-2})-2\det(D_{n-3})-4\det(D_{n-4}))]-2\det(D_{n-1})-4\det(D_{n-2})\\ 
	& =4\det(D_{n-2})-16\det ({{D}_{n-3}})-32\det ({{D}_{n-4}}), 
\end{aligned}
\end{equation*}
i.e.,
\begin{equation}\label{eqn-rr-S5-T52n-1}
    \det ({{D}_{n}})-4\det(D_{n-2})+16\det ({{D}_{n-3}})+32\det ({{D}_{n-4}})=0.
\end{equation}

Therefore, we obtain a third five-term recurrence relation for the determinant of the distance matrix of trees.


\subsection{The recurrence relation related to $T_n=T_{5,2}\circ_{n} R$}\label{subsec-S5-T52n}

In this subsection, we consider the decomposition $T_n=T_{5,2}\circ_{n} R$ of $T_n$, and we interchange the labels of $n-1$ and $n-3$ and interchange the labels of $n-4$ and $n$, then the partitioned matrix $D_n$ has the following form:
\begin{equation}\label{eqn-S5-T52n}
 {{D}_{n}} =\left( \begin{matrix}
		{{D}_{n-5}} & \alpha  & \alpha +e & \alpha +2e & \alpha +3e & \alpha +3e  \\
		{{\alpha }^{T}} & 0 & 1 & 2 & 3 & 3  \\
		{{\alpha }^{T}}+{{e}^{T}} & 1 & 0 & 1 & 2 & 2  \\
		{{\alpha }^{T}}+2{{e}^{T}} & 2 & 1 & 0 & 1 & 1  \\
		{{\alpha }^{T}}+3{{e}^{T}} & 3 & 2 & 1 & 0 & 2  \\
		{{\alpha }^{T}}+3{{e}^{T}} & 3 & 2 & 1 & 2 & 0  \\
\end{matrix} \right),
\end{equation}
where $\alpha$ is a column vector whose entries are $d_{i,n-4}$ for $i\in [n-5]$.

We take ${{R}_{3}}-2{{R}_{4}}\rightarrow {{R}_{6}}$ and expand it by the last row, take $\frac{1}{2} ({{C}_{3}}- {{C}_{5}})\rightarrow {{C}_{4}}$ on the first determinant of the second line, expand the second determinant of the fourth line, take ${{C}_{2}}-{{C}_{3}}\rightarrow {{C}_{4}}$ on the two determinants of the fifth line and take ${{R}_{2}}-{{R}_{3}}\rightarrow {{R}_{4}}$ on the second determinant of the same line, then we have 
\begin{equation*}
\begin{aligned}
&\hspace{-6em}	 \det \left( {{D}_{n}} \right)=\left| \begin{matrix}
		{{D}_{n-5}} & \alpha  & \alpha +e & \alpha +2e & \alpha +3e & \alpha +3e  \\
		{{\alpha }^{T}} & 0 & 1 & 2 & 3 & 3  \\
		{{\alpha }^{T}}+{{e}^{T}} & 1 & 0 & 1 & 2 & 2  \\
		{{\alpha }^{T}}+2{{e}^{T}} & 2 & 1 & 0 & 1 & 1  \\
		{{\alpha }^{T}}+3{{e}^{T}} & 3 & 2 & 1 & 0 & 2  \\
		0 & 0 & 0 & 2 & 2 & 0  \\
	\end{matrix} \right| \\ 
&\hspace{-6em}	=2\left| \begin{matrix}
		{{D}_{n-5}} & \alpha  & \alpha +e & \alpha +3e & \alpha +3e  \\
		{{\alpha }^{T}} & 0 & 1 & 3 & 3  \\
		{{\alpha }^{T}}+{{e}^{T}} & 1 & 0 & 2 & 2  \\
		{{\alpha }^{T}}+2{{e}^{T}} & 2 & 1 & 1 & 1  \\
		{{\alpha }^{T}}+3{{e}^{T}} & 3 & 2 & 0 & 2  \\
	\end{matrix} \right|-2\left| \begin{matrix}
		{{D}_{n-5}} & \alpha  & \alpha +e & \alpha +2e & \alpha +3e  \\
		{{\alpha }^{T}} & 0 & 1 & 2 & 3  \\
		{{\alpha }^{T}}+{{e}^{T}} & 1 & 0 & 1 & 2  \\
		{{\alpha }^{T}}+2{{e}^{T}} & 2 & 1 & 0 & 1  \\
		{{\alpha }^{T}}+3{{e}^{T}} & 3 & 2 & 1 & 2  \\
	\end{matrix} \right|  \\
&\hspace{-6em}  =2\left| \begin{matrix}
		{{D}_{n-5}} & \alpha  & \alpha +e & \alpha +2e & \alpha +3e  \\
		{{\alpha }^{T}} & 0 & 1 & 2 & 3  \\
		{{\alpha }^{T}}+{{e}^{T}} & 1 & 0 & 1 & 2  \\
		{{\alpha }^{T}}+2{{e}^{T}} & 2 & 1 & 1 & 1  \\
		{{\alpha }^{T}}+3{{e}^{T}} & 3 & 2 & 0 & 2  \\
	\end{matrix} \right|-2\det ({{D}_{n-1}})-4\det ({{D}_{n-2}}) \\ 
	&\hspace{-6em} =2\left| \begin{matrix}
		{{D}_{n-5}} & \alpha  & \alpha +e & \alpha +2e & \alpha +3e  \\
		{{\alpha }^{T}} & 0 & 1 & 2 & 3  \\
		{{\alpha }^{T}}+{{e}^{T}} & 1 & 0 & 1 & 2  \\
		{{\alpha }^{T}}+2{{e}^{T}} & 2 & 1 & 0 & 1  \\
		{{\alpha }^{T}}+3{{e}^{T}} & 3 & 2 & 1 & 2  \\
	\end{matrix} \right|+2\left| \begin{matrix}
		{{D}_{n-5}} & \alpha  & \alpha +e & 0 & \alpha +3e  \\
		{{\alpha }^{T}} & 0 & 1 & 0 & 3  \\
		{{\alpha }^{T}}+{{e}^{T}} & 1 & 0 & 0 & 2  \\
		{{\alpha }^{T}}+2{{e}^{T}} & 2 & 1 & 1 & 1  \\
		{{\alpha }^{T}}+3{{e}^{T}} & 3 & 2 & -1 & 2  \\
	\end{matrix} \right| -2\det ({{D}_{n-1}})-4\det ({{D}_{n-2}}) \\ 
	&\hspace{-6em} = 2\left| \begin{matrix}
		{{D}_{n-5}} & \alpha  & \alpha +e & \alpha +3e  \\
		{{\alpha }^{T}} & 0 & 1 & 3  \\
		{{\alpha }^{T}}+{{e}^{T}} & 1 & 0 & 2  \\
		{{\alpha }^{T}}+2{{e}^{T}} & 2 & 1 & 1  \\
	\end{matrix} \right|+2\left| \begin{matrix}
		{{D}_{n-5}} & \alpha  & \alpha +e & \alpha +3e  \\
		{{\alpha }^{T}} & 0 & 1 & 3  \\
		{{\alpha }^{T}}+{{e}^{T}} & 1 & 0 & 2  \\
		{{\alpha }^{T}}+3{{e}^{T}} & 3 & 2 & 2  \\
	\end{matrix} \right| \\
&\hspace{-6em}  =2\left| \begin{matrix}
		{{D}_{n-5}} & \alpha  & \alpha +e & \alpha +2e  \\
		{{\alpha }^{T}} & 0 & 1 & 2  \\
		{{\alpha }^{T}}+{{e}^{T}} & 1 & 0 & 3  \\
		{{\alpha }^{T}}+2{{e}^{T}} & 2 & 1 & 2  \\
	\end{matrix} \right|+2\left| \begin{matrix}
		{{D}_{n-5}} & \alpha  & \alpha +e & \alpha +2e  \\
		{{\alpha }^{T}} & 0 & 1 & 2  \\
		{{\alpha }^{T}}+{{e}^{T}} & 1 & 0 & 3  \\
		{{\alpha }^{T}}+2{{e}^{T}} & 2 & 3 & 2  \\
	\end{matrix} \right| 
\end{aligned}
\end{equation*}
\begin{equation*}
\begin{aligned}
&  =4\left| \begin{matrix}
		{{D}_{n-5}} & \alpha  & \alpha +e & \alpha +2e  \\
		{{\alpha }^{T}} & 0 & 1 & 2  \\
		{{\alpha }^{T}}+{{e}^{T}} & 1 & 0 & 3  \\
		{{\alpha }^{T}}+2{{e}^{T}} & 2 & 1 & 2  \\
	\end{matrix} \right|+2\left| \begin{matrix}
		{{D}_{n-5}} & \alpha  & \alpha +e & \alpha +2e  \\
		{{\alpha }^{T}} & 0 & 1 & 2  \\
		{{\alpha }^{T}}+{{e}^{T}} & 1 & 0 & 3  \\
		0 & 0 & 2 & 0  \\
	\end{matrix} \right|\\
&=4\det ({{D}_{n-2}})+4\left| \begin{matrix}
		{{D}_{n-5}} & \alpha  & \alpha +e & 0  \\
		{{\alpha }^{T}} & 0 & 1 & 0  \\
		{{\alpha }^{T}}+{{e}^{T}} & 1 & 0 & 2  \\
		{{\alpha }^{T}}+2{{e}^{T}} & 2 & 1 & 2  \\
	\end{matrix} \right|-4\left| \begin{matrix}
		{{D}_{n-5}} & \alpha  & \alpha +2e  \\
		{{\alpha }^{T}} & 0 & 2  \\
		{{\alpha }^{T}}+{{e}^{T}} & 1 & 3  \\
	\end{matrix} \right| \\ 
	& =4\det ({{D}_{n-2}})+8\det ({{D}_{n-3}})-8\left| \begin{matrix}
		{{D}_{n-5}} & \alpha  & \alpha +e  \\
		{{\alpha }^{T}} & 0 & 1  \\
		{{\alpha }^{T}}+2{{e}^{T}} & 2 & 1  \\
	\end{matrix} \right|-4\left| \begin{matrix}
		{{D}_{n-5}} & \alpha  & \alpha +2e  \\
		{{\alpha }^{T}} & 0 & 2  \\
		{{\alpha }^{T}}+{{e}^{T}} & 1 & 3  \\
	\end{matrix} \right| \\ 
	& =4\det ({{D}_{n-2}})+8\det ({{D}_{n-3}})-16\left| \begin{matrix}
		{{D}_{n-5}} & \alpha  & \alpha +e  \\
		{{\alpha }^{T}} & 0 & 1  \\
		{{\alpha }^{T}}+{{e}^{T}} & 1 & 1  \\
	\end{matrix} \right|-8\left| \begin{matrix}
		{{D}_{n-5}} & \alpha  & \alpha +2e  \\
		{{\alpha }^{T}} & 0 & 1  \\
		{{\alpha }^{T}}+{{e}^{T}} & 1 & 2  \\
	\end{matrix} \right| \\ 
	& =4\det ({{D}_{n-2}})+8\det ({{D}_{n-3}})-16[\det(D_{n-3})+\det(D_{n-4})]-8[\det(D_{n-3})+2\det(D_{n-4})]\\ 
	& =4\det ({{D}_{n-2}})-16\det ({{D}_{n-3}})-32\det ({{D}_{n-4}})  
\end{aligned}
\end{equation*}
i.e.,
\begin{equation}\label{eqn-rr-S5-T52n}
\det ({{D}_{n}})-4\det ({{D}_{n-2}})+16\det ({{D}_{n-3}})+32\det ({{D}_{n-4}})=0.
\end{equation}

This recurrence relation is the same as Equation~(\ref{eqn-rr-S5-T52n-1}).


\subsection{The recurrence relation related to $T_n=K_{1,4}\circ_{n-4} R$}\label{subsec-S5-K14n-4}

In this subsection, we consider the decomposition $T_n=K_{1,4}\circ_{n-4} R$ of $T_n$, then the partitioned matrix $D_n$ has the following form:
\begin{equation}\label{eqn-S5-K14n-4}
 {{D}_{n}} =\left( \begin{matrix}
		{{D}_{n-5}} & \alpha  & \alpha +e & \alpha +e & \alpha +e & \alpha +e  \\
		{{\alpha }^{T}} & 0 & 1 & 1 & 1 & 1  \\
		{{\alpha }^{T}}+{{e}^{T}} & 1 & 0 & 2 & 2 & 2  \\
		{{\alpha }^{T}}+{{e}^{T}} & 1 & 2 & 0 & 2 & 2  \\
		{{\alpha }^{T}}+{{e}^{T}} & 1 & 2 & 2 & 0 & 2  \\
		{{\alpha }^{T}}+{{e}^{T}} & 1 & 2 & 2 & 2 & 0  \\
\end{matrix} \right), 
\end{equation}
where $\alpha$ is a column vector whose entries are $d_{i,n-4}$ for $i\in [n-5]$. 

We take $-{{R}_{5}}\rightarrow {{R}_{6}}$ and expand it by the last row, then we get 
\begin{equation*}
\begin{aligned}
	 \det \left( {{D}_{n}} \right)&=\left| \begin{matrix}
		{{D}_{n-5}} & \alpha  & \alpha +e & \alpha +e & \alpha +e & \alpha +e  \\
		{{\alpha }^{T}} & 0 & 1 & 1 & 1 & 1  \\
		{{\alpha }^{T}}+{{e}^{T}} & 1 & 0 & 2 & 2 & 2  \\
		{{\alpha }^{T}}+{{e}^{T}} & 1 & 2 & 0 & 2 & 2  \\
		{{\alpha }^{T}}+{{e}^{T}} & 1 & 2 & 2 & 0 & 2  \\
		0 & 0 & 0 & 0 & 2 & -2  \\
	\end{matrix} \right| \\ 
	& =-2\det ({{D}_{n-1}})-2[\det ({{D}_{n-1}})+2\det ({{D}_{n-2}})] \\ 
	& =-4\det ({{D}_{n-1}})-4\det ({{D}_{n-2}})  
\end{aligned}
\end{equation*}
i.e.,
\begin{equation}\label{eqn-rr-S5-K14n-4}
	\det \left( {{D}_{n}} \right)+4\det \left( {{D}_{n-1}} \right)+4\det \left( {{D}_{n-2}} \right)=0.
\end{equation}


\subsection{The recurrence relation related to $T_n=K_{1,4}\circ_{n-3} R$}\label{subsec-S5-K14n-3}

In this subsection, we consider the decomposition $T_n=K_{1,4}\circ_{n-3} R$ of $T_n$, and we interchange the labels of $n-4$ and $n-3$, then the partitioned matrix $D_n$ has the following form:
\begin{equation}\label{eqn-S5-K14n-3}
 {{D}_{n}} =\left( \begin{matrix}
	{{D}_{n-5}} & \alpha  & \alpha +e & \alpha +2e & \alpha +2e & \alpha +2e  \\
	{{\alpha }^{T}} & 0 & 1 & 2 & 2 & 2  \\
	{{\alpha }^{T}}+{{e}^{T}} & 1 & 0 & 1 & 1 & 1  \\
	{{\alpha }^{T}}+2{{e}^{T}} & 2 & 1 & 0 & 2 & 2  \\
	{{\alpha }^{T}}+2{{e}^{T}} & 2 & 1 & 2 & 0 & 2  \\
	{{\alpha }^{T}}+2{{e}^{T}} & 2 & 1 & 2 & 2 & 0  \\
\end{matrix} \right), 
\end{equation}
where $\alpha$ is a column vector whose entries are $d_{i,n-4}$ for $i\in [n-5]$. When we swap $R_2$ and $R_3$ and swap $C_2$ and $C_3$, then we get the same matrix as in Equation~(\ref{eqn-S5-K14n-4}). Hence we can get the same recurrence relation as in Equation~(\ref{eqn-rr-S5-K14n-4}).

%


\subsection{Summary of recurrence relations related to $S_5$}\label{subsec-S5-summary}

In this section, we considered all the nine cases of decompositions $T_n=S_5\circ R$, and we get the recurrence relations listed in Table~\ref{tab-S5}.

\begin{table}[ht]
	\caption{Recurrence relations related to $T_n=S_5\circ R$}
	\resizebox{\columnwidth}{!}{
	\label{tab-S5}
	\begin{tabular}{|c|c|c|}
		\hline
	decomposition & recurrence relation & remark\\ \hline
	\begin{tabular}{c}$T_n=P_5\circ_{n-4} R$, or\\ $T_n=P_5\circ_{n-3} R$~~~~\end{tabular} & $\det ({{D}_{n}})+2\det ({{D}_{n-1}})+8\det ({{D}_{n-3}})+16\det ({{D}_{n-4}})=0$  & five-term\\ \hline
	$T_n=P_5\circ_{n-2} R$	&   $\det ({{D}_{n}})+2\det ({{D}_{n-1}})-2\det ({{D}_{n-2}})+8\det ({{D}_{n-4}})=0$ & five-term  \\ \hline
	\begin{tabular}{c}$T_n=T_{5,2}\circ_{n-4} R$, or\\ $T_n=T_{5,2}\circ_{n-2} R$~~~~\end{tabular}	&   $\det ({{D}_{n}})+2\det ({{D}_{n-1}})-4\det ({{D}_{n-2}})-8\det ({{D}_{n-3}})=0$ & four-term \\ \hline 
	\begin{tabular}{c}$T_n=T_{5,2}\circ_{n-1} R$, or\\ $T_n=T_{5,2}\circ_{n} R$~~~~\end{tabular}	&   $\det ({{D}_{n}})-4\det ({{D}_{n-2}})+16\det ({{D}_{n-3}})+32\det ({{D}_{n-4}})=0$ & five-term \\ \hline 
	\begin{tabular}{c}$T_n=K_{1,4}\circ_{n-4} R$, or\\ $T_n=K_{1,4}\circ_{n-3} R$~~~~\end{tabular}	&   $\det({{D}_{n}})+4\det({{D}_{n-1}})+4\det({{D}_{n-2}})=0$ & three-term \\ \hline 
	\end{tabular}}
\end{table}


\section{Recurrence relations related to $T_n=S_6\circ R$}\label{sec-S6}

In this section, we consider the decomposition $T_n=S_6\circ R$ and recurrence relations related. We may suppose the vertex set of the subtree $S_6$ on five vertices is $[n]\setminus [n-6]=\{n-5,n-4,n-3,n-2,n-1,n\}$. We know the subtree $S_6$ is a path $P_6$, star-like trees $K(1,1,3)$, $K(1,2,2)$, $K(1,1,1,2)$, a dumbbell tree $T_{6,5}$ or a star $K_{1,5}$, where $P_6$ has edges $\{n-5,n-4\}$, $\{n-4,n-3\}$, $\{n-3,n-2\}$, $\{n-2,n-1\}$ and $\{n-1,n\}$, $K(1,1,3)$ has edges $\{n-5,n-3\}$, $\{n-4,n-3\}$, $\{n-3,n-2\}$, $\{n-2,n-1\}$ and $\{n-1,n\}$, $T_{6,5}$ has edges $\{n-5,n-3\}$, $\{n-4,n-3\}$, $\{n-3,n-2\}$, $\{n-2,n-1\}$ and $\{n-2,n\}$. 

There are three different decompositions with form $T_n=P_6\circ R$, five different decompositions with form $T_n=K(1,1,3)\circ R$, four different decompositions with form $T_n=K(1,2,2)\circ R$, four different decompositions with form $T_n=K(1,1,1,2)\circ R$, two different decompositions with form $T_n=T_{6,5}\circ R$, two different decompositions with form $T_n=K_{1,5}\circ R$. Hence there are twenty different decompositions of the form $T_n=S_6\circ R$. 

As for the length of this manuscript, we omit the details of computations for the recurrence relations on the determinants of trees of these twenty cases, as a result the following listed recurrence relations in Table~\ref{tab-S6} were obtained.

\begin{table}[ht]
	\caption{Recurrence relations related to $T_n=S_6\circ R$}
	\resizebox{\columnwidth}{!}{
	\label{tab-S6}
	\begin{tabular}{|c|c|c|}
%
	\hline
	decomposition  & recurrence relation &  remark  \\ \hline
	$P_6\circ_{n-5} R$	& $\det \left( {{D}_{n}} \right)+2\det \left( {{D}_{n-1}} \right)-16\det \left( {{D}_{n-4}} \right)-32\det \left( {{D}_{n-5}} \right)=0$ & six-term \\ \hline
	$P_6\circ_{n-3} R$	& $\det ({{D}_{n}})+2\det ({{D}_{n-1}})-4\det ({{D}_{n-2}})-16\det ({{D}_{n-3}})-32\det ({{D}_{n-4}})-32\det ({{D}_{n-5}})=0$ & six-term  \\ \hline
	$K(1,1,3)\circ_{n-1} R$	& $\det ({{D}_{n}})-2\det ({{D}_{n-1}})-4\det ({{D}_{n-2}})+16\det ({{D}_{n-3}})-32\det ({{D}_{n-4}})-96\det ({{D}_{n-5}})=0$ & six-term \\ \hline
	$K(1,1,3)\circ_{n} R$	&  $\det ({{D}_{n}})-2\det ({{D}_{n-1}})-8\det ({{D}_{n-2}})+24\det ({{D}_{n-3}})+48\det ({{D}_{n-4}})=0$ & five-term  \\ \hline

	$2$ decompositions	& $\det ({{D}_{n}})+2\det ({{D}_{n-1}})+8\det ({{D}_{n-3}})+16\det ({{D}_{n-4}})=0$ & five-term  \\ \hline
	$7$ decompositions	& $\det ({{D}_{n}})+2\det ({{D}_{n-1}})-4\det ({{D}_{n-2}})-8\det ({{D}_{n-3}})=0$ & four-term  \\ \hline
	$7$ decompositions	& $\det ({{D}_{n}})+4\det ({{D}_{n-1}})+4\det ({{D}_{n-2}})=0$  & three-term  \\ \hline
	\end{tabular}}
\end{table}


\section{Recurrence relations related to $T_n=S_7\circ R$}\label{sec-S7}

In this section, we consider the decomposition $T_n=S_7\circ R$ and recurrence relations related. We may suppose the vertex set of the subtree $S_6$ on five vertices is $[n]\setminus [n-7]=\{n-6,n-5,n-4,n-3,n-2,n-1,n\}$. We know the subtree $S_7$ is a path $P_7$, star-like trees $K(1,1,4)$, $K(1,2,3)$, $K(2,2,2)$, $K(1,1,1,3)$, $K(1,1,2,2)$, $K(1,1,1,2,2)$, three trees $T_{6,5}^{1,1}$, $T_{6,5}^{1,2}$, $T_{6,5}^2$ obtained from the dumbbell tree $T_{6,5}$ by adding a new vertex with degree one to a pendent vertex, adding a new vertex with degree one to a vertex with degree two, or adding a new vertex between two vertices with degree three, respectively, or a star $K_{1,6}$, where $P_6$ has edges $\{n-6,n-5\}$, $\{n-5,n-4\}$, $\{n-4,n-3\}$, $\{n-3,n-2\}$, $\{n-2,n-1\}$ and $\{n-1,n\}$. 

There are four different decompositions with form $T_n=P_7\circ R$, six different decompositions with form $T_n=K(1,1,4)\circ R$, seven different decompositions with form $T_n=K(1,2,3)\circ R$, three different decompositions with form $T_n=K(2,2,2)\circ R$, five different decompositions with form $T_n=K(1,1,1,3)\circ R$, four different decompositions with form $T_n=K(1,1,2,2)\circ R$, four different decompositions with form $T_n=K(1,1,1,1,2)\circ R$, six different decompositions with form $T_n=T_{6,5}^{1,1}\circ R$, four different decompositions with form $T_n=T_{6,5}^{1,2}\circ R$, three different decompositions with form $T_n=T_{6,5}^2\circ R$, two different decompositions with form $T_n=K_{1,6}\circ R$. Hence there are forty-eight different decompositions of the form $T_n=S_6\circ R$. 

As for the length of this manuscript, we omit the details of computations for the recurrence relations on the determinants of trees of these forty-eight cases, as a result the following listed recurrence relations in Table~\ref{tab-S7} were obtained.

\begin{table}[ht]
	\caption{Recurrence relations related to $T_n=S_7\circ R$}
	\resizebox{\columnwidth}{!}{
	\label{tab-S7}
	\begin{tabular}{|c|c|c|}
	\hline
   decomposition & recurrence relation & remark \\ \hline
	\begin{tabular}{c}$T_n=P_7\circ_{n-6} R$, or\\ $T_n=P_7\circ_{n-5} R$~~~~\end{tabular}	&  $\det \left( {{D}_{n}} \right)+2\det \left( {{D}_{n-1}} \right)+32\det \left( {{D}_{n-5}} \right)+64\det \left( {{D}_{n-6}} \right)=0$  & seven-term \\ \hline 
	$T_n=P_7\circ_{n-4} R$	&  $\det \left( {{D}_{n}} \right)+2\det \left( {{D}_{n-1}} \right)+8\det \left( {{D}_{n-3}} \right)+64\det \left( {{D}_{n-4}} \right)+192\det \left( {{D}_{n-5}} \right)+192\det \left( {{D}_{n-6}} \right)=0$   & seven-term  \\ \hline
	$T_n=P_7\circ_{n-3} R$	&  $\det \left( {{D}_{n}} \right)+2\det \left( {{D}_{n-1}} \right)+8\det \left( {{D}_{n-3}} \right)+32\det \left( {{D}_{n-4}} \right)+64\det \left( {{D}_{n-5}} \right)+64\det \left( {{D}_{n-6}} \right)=0$   & seven-term  \\ \hline
	$2$ decompositions	&  $\det \left( {{D}_{n}} \right)+2\det \left( {{D}_{n-1}} \right)-16\det \left( {{D}_{n-4}} \right)-32\det \left( {{D}_{n-5}} \right)=0$    & six-term \\ \hline
	$4$ decompositions	&  $\det ({{D}_{n}})+2\det ({{D}_{n-1}})+8\det ({{D}_{n-3}})+16\det ({{D}_{n-4}})=0$    & five-term \\ \hline
	$18$ decompositions	&  $\det ({{D}_{n}})+2\det ({{D}_{n-1}})-4\det ({{D}_{n-2}})-8\det ({{D}_{n-3}})=0$    & four-term \\ \hline
	$20$ decompositions	&  $\det ({{D}_{n}})+4\det ({{D}_{n-1}})+4\det ({{D}_{n-2}})=0$    & three-term \\ \hline
	\end{tabular}}
\end{table}

\end{document}